\documentclass[12pt,leqno]{article}
\usepackage{amsmath,amssymb}
\usepackage{amsfonts}
\usepackage{amsthm}
\usepackage{newlfont}
\usepackage{amssymb}
\usepackage{amsthm}
\usepackage{newlfont}
\usepackage{mathrsfs}
\usepackage{amsmath,hhline}
\usepackage{setspace}
\usepackage [latin1]{inputenc}
\usepackage{graphicx}

\newtheorem{theorem}{Theorem}[section]

\newtheorem{remark}[theorem]{Remark}

\newtheorem{lemma}[theorem]{Lemma}
\newtheorem{claim}{Claim}

\newcommand{\RR}{{\mathbb R}}

\begin{document}
\singlespacing

\title{Radial and nonradial solutions for nonlinear scalar field equations and a new proof for certain Gagliardo-Nirenberg type inequality on the plane}

\author{{\Large Qiuping Lu}\footnote{School of Mathematical Science, Yangzhou University, Yangzhou 225002, China.
 {\tt qplu@yzu.edu.cn} \ Supported
by P.R.C. grant NSFC 11301460 and NSFC 12171418.} \and  {\Large Shuolei Liu\footnote{School of Mathematical Science, Yangzhou University, Yangzhou 225002, China. {\tt 13781550416@163.com} 
\ Supported by P.R.C. grant BK 20191436 and BK 20221369 .}
 }}


\date{ }

\maketitle

\pagenumbering{arabic}
\begin{flushleft}
 \textbf{Abstract.} A nonlinear scalar field equation is studied, which takes the form
 \begin{equation*}
 - \triangle{u}=g(u) \ \textrm{in} \ \RR^{N}, \quad \quad u \in E,
 \end{equation*}
 where $E$ is defined in \eqref{main}. The nonlinear term $g(s)$ has subcritical growth at infinity. Near $s=0$, the nonlinear term $g(s)$ grows like $s^q$ for any $q>0$ if $N=2$ and for any $q \in (0,\ (N+2)/(N-2))$ if $N\ge 3$.
 Through the technique of adding one dimension of space, by mountain pass theorem and symmetric mountain theorem, we prove the multiplicity of radial solutions if $N\ge 2$ and the multiplicity of nonradial solutions if $N\ge 4$. Unlike \cite{MED} and \cite{JLU}, the proof presented here does not rely on any decomposition results. During the treatment for $N=2$, we provide a new proof for certain Gagliardo-Nirenberg type inequality on $\RR^2$, which is originally due to \cite{EI}. Based on this inequality we get a Trudinger-Moser type inequality involving gradient norm through Taylor expansion. Eventually we construct a special path to obtain the least energy level for solutions.


\textbf{Data Availability Statement.} This is a purely mathematical research article and data sharing is not applicable to this article as no new data were created or analyzed in this study.
\end{flushleft}
\section{Introduction}
In this paper, for $N \ge 2$, we want to find radial and nonradial solutions of the following nonlinear scalar field equation,
\begin{equation}\label{main}
- \triangle{u}=g(u) \ \textrm{in} \ \RR^{N}, \quad \quad u \in E,
\end{equation}
where $E$ is the completion of $ C_0^{\infty}(\RR^{N})$ under the norm
\begin{equation*}
||v||_{E} := (\int_{\RR^{N}} |\nabla v|^{2} dx)^{\frac{1}{2}} + (\int_{\RR^{N}} |v|^{q+1} dx)^{\frac{1}{q+1}}, \quad \textrm{ for any $v \in C_0^{\infty}(\RR^{N}) $},
\end{equation*}
and $g : \RR \to \RR$ is a continuous and odd function, which satisfies the following conditions:
\begin{enumerate}
\item [(g1)] $ -\infty < \liminf_{s \to 0^{+}} g(s)/ s^q \le \limsup_{s \to 0^{+}} g(s)/ s^q =-2m<0$ where $m>0$, $0<q< 2^{*}-1 $, $2^{*}=\infty $ if $N=2$ and $2^*=(2N)/(N-2)$ if $N\ge 3$, and the factor $2$ in front of $m$ is just for convenience;
\item [(g2)] if $N=2$,\ for any $\alpha>0$, there exists $C_{\alpha}>0$ such that $ |g(s)|\le C_{\alpha}e^{\alpha s^2} \ \textrm{for all $s\in \RR$} $,
             if $N\ge 3$, $ \lim_{s\to \infty}|g(s)| / |s|^{2^{*}-1}=0$;
\item [(g3)] $G(s_{0})>0$ for some $s_{0}>0$, where $G(s):= \int_{0}^{s}g(t)dt $.
\end{enumerate}

This type of nonlinear problem arises in many models in mathematical physics. For example \eqref{main} is encountered during the search for certain kinds of solitary waves
in nonlinear Klein-Gordon equations or Schr\"{o}dinger equations, we refer to \cite{BL1} for more details.
When $q=1$ and $N\ge 3$, under the assumptions (g1), (g2) and (g3), Berestycki and Lions \cite{BL1} proved the existence of a radial least energy solution
by showing that the following constrained minimization problem
\begin{equation}\label{constraint1}
\min  \{  \int_{\RR^N}  |\nabla{v}|^2 dx \ \ \ \ | \ \ \ \int_{\RR^N}  G(v)dx =1 \ \ \textrm{and}  \ \ v \in E \}
\end{equation}
had a radial solution, and they \cite{BL2} also studied the existence of infinitely many radial solutions through minimax methods, see also \cite{ST} for related work.
When $q=1$ and $N=2$, under (g2), (g3) and the more restrictive assumption $-\infty<\lim_{s\to 0} g(s)/s <0$, Berestycki, Gallou\"{e}t and Kavian \cite{BGK} obtained the existence of a radial least energy solution
by proving that the following constrained minimization problem
\begin{equation}\label{constraint2}
\min  \{  \int_{\RR^2}  |\nabla{v}|^2 dx \ \ \ \ | \ \ \ \int_{\RR^2}  G(v)dx =0 \ \ \textrm{and}  \ \ v \in E-\{0\}  \}
\end{equation}
admitted a radial solution, and the existence of infinitely many solutions were also studied. We point out that the difference between \eqref{constraint1} and \eqref{constraint2} is mainly due to
the Pohozaev identity satisfied by the solutions of \eqref{main}. See also Strauss \cite{STR} and Coleman, Glaser and Martin \cite{CGM} for earlier works when $q=1$ and $N\ge 2$. When $N\ge3$ and $0<q<1$,
with (g1), (g2) and (g3), Lu and Ling \cite{SQ5} proved the existence of a radial least energy solution to \eqref{main} with compact support,
and they used this solution to construct least energy solution to certain singularly perturbed problem in a bounded domain and got totally different peaking phenomenon from those in Ni and Takagi \cite{NT1} and \cite{NT2},
Ni and Wei \cite{NW1}. For the study of related system equations, we refer to Brezis and Lieb \cite{BLI}  and Bruning \cite{BRU}.

When it comes to the question of the existence of nonradial solutions, which was raised by Berestycki and Lions \cite{BL1} at the end of section $10.8$, Bartsch and Willem \cite{BW} proved the existence of nonradial solution if $q=1$, $N\ge 4$ and $N\neq 5$ under Ambrosetti-Rabinowitz type condition \cite{AR}. Lorca and Ubilla \cite{LU} figured out the existence of nonradial solution if $q=1$ and $N=5$. Under (g1), (g2) and (g3), Mederski \cite{MED} developed a new critical point theory for constrained functionals and showed the existence of nonradial solutions without Ambrosetti-Rabinowitz type condition if $q=1$ and $N\ge 4$. Afterwards, under the same assumptions, Jeanjean and Lu \cite{JLU} recovered the same results by using monotonicity trick. Finally Musso, Pacard and Wei \cite{MPW} obtained noradial solutions for any $N\ge 2$, see also \cite{AMP}. But in their works, certain nondegeneracy condition is needed to apply Lyapunov-Schmidt type reduction process and this nondegeneracy condition is much stronger than (g1), (g2) and (g3).

We define the energy functional associated to \eqref{main} by
\begin{equation*}
I(v):= \frac{1}{2} \int_{\RR^{N}} |\nabla{v}|^{2}dx- \int_{\RR^{N}}G(v)dx, \quad v \in E.
\end{equation*}For $q=1$ and $N\ge 2$, assuming (g1), (g2) and (g3), Jeanjean and Tanaka \cite{JT} found that $v=0$ was a strict local minimum point of the energy functional $I(v)$
and $I(v)$ had a mountain pass structure around $v=0$. Based on this fact, by the technique of adding one dimension of space \cite{JEA}, Hirata, Ikoma and Tanaka \cite{HIT} showed the existence of infinitely many radial solutions including a least energy solution through mountain pass and symmetric mountain approaches, therefore \cite{HIT} is an improvement of the existence results in \cite{BGK} by assuming (g1) instead of $-\infty<\lim_{s\to 0} g(s)/s <0$.

In this paper, relying on the technique of adding one dimension of space, by working on different Banach spaces we are going to recover the same results in \cite{HIT} and \cite{MED}, and extend these existence results to the case $0<q< 2^{*}-1$ for $N\ge 2$. Unlike \cite{MED} and \cite{JLU}, we emphasis that our proof will not use any decomposition results.

In order to extend these results, we need more estimates if $N=2$. Indeed if $N\ge 3$, according to standard Sobolev embedding theorems, (g2) implies that $g(s)$ has subcritical growth at $\infty$. If $N=2$, we also say $g$ has subcritical growth at $\infty$ if $g$ satisfies
\begin{enumerate}
\item [(g2e)] \ \ \ \ \ \ \ \ \ \ \ \ \ \ \ \ \ \ \ \ \ \ \ \ $ \lim_{s\to \infty}\frac{|g(s)|}{e^{\alpha s^{2}}}=0$ for any $\alpha>0$.
\end{enumerate}
The classification of criticality is motivated by the Trudinger-Moser inequality \cite{TRU,MOS}, which basically says that for a bounded domain $\Omega\subset \RR^2$, $\int_{\Omega} e^{v^2} dx $ is finite if $v\in H^{1}_{0}(\Omega)$. It is clear that (g2) is equivalent to (g2e) if $N=2$, so overall (g2) means that the nonlinear term $g$ of \eqref{main} has subcritical growth at $\infty$.

Now We present a sharp Gagliardo-Nirenberg type inequality due to Edmunds and Ilyin \cite{EI}, and the coefficient $C(q,p)$ below is given in terms of gamma function, which is equivalent to the one in \cite{EI}. Based on this inequality and Taylor expansion formula, we derive an interesting Trudinger-Moser type inequality. This Trudinger-Moser type inequality is not in optimal form and we will deal with this problem in a forthcoming paper. For simplicity, we write $(\int_{\RR^{N}} |v|^{p} dx)^{\frac{1}{p}}$ as $ ||v||_{p}$ when $p\ge 1$.
\begin{theorem}\label{TRU-MOS}
(i) \ Let $v \in C_{0}^{\infty}(\RR^2) $ and $p \in [ q+1, \ \infty)$.  Then
$$ ||v||_{p} \le C(q,p) ||\nabla{v}||_{2}^{1-\frac{q+1}{p}} ||v||_{q+1}^{\frac{q+1}{p}} ,$$
where the coefficient $C(q,p)=2(1+\Gamma(p/2 +1))^{1/p}$, $C(q,p)\le Q p^{1/2}$ and $Q$ is a constant depending only on $q$.

(ii)\ There exists a small constant $\rho_1>0$ depending only on $q$ such that, for any nonnegative $k$ with $k\ge q-1$,
 $$ \int_{\RR^2}|v|^{k}(e^{v^2} -1)dx \le C ||\nabla{v}||_{2}^{k+1-q}||v||_{q+1}^{q+1}  \ \ \ \   \textrm{for any $v\in E$ with  $||\nabla{v}||_2 \le \rho_1 $ ,}$$
where $C$ depends on $k$ and $q$.
\end{theorem}
Classical Trudinger-Moser inequality is over bounded domain, the Trudinger-Moser type inequality we developed here is over $\RR^2$, there are lots of works in this direction,
interested readers may check Cao \cite{CAO}, Do \'O \cite{JMB}, Ozawa \cite{OZA}, Adachi and Tanaka \cite{AT} and Ruf \cite{RUF},
we refer to Sandeep \cite{SAN} for a more complete set of references about Trudinger-Moser type inequality over bounded and unbounded domains.

By (g1), (g2), standard Sobolev embedding inequality for $N\ge 3 $ and theorem~\ref{TRU-MOS}, for any $v\in E$, $\int_{\RR^{N}}G(v)dx$ is finite, therefore the energy functional $I(v)$ does make sense, see lemma~\ref{makesense} for details.
Moreover theorem~\ref{c1functional} shows that $I:E \to \RR$ is of $C^{1}$.
Now searching for solutions of \eqref{main} becomes finding critical points of $I(v)$, by lemm~\ref{mountain-pass-structure-1}, the energy functional $I(v)$ has a mountain pass structure around $v=0$. So we can define a sequence of minimax value $c_n^H$:
\begin{equation*}
c_n^H:=\inf_{\gamma \in \Gamma_{n}^H} \sup_{x\in D_n} I(\gamma(x)),\ \ \ D_n :=\{ x=(x_1,\ x_2,\ \cdots,\ x_n)  \in \RR^n \ |\ |x|^2=\Sigma_{i=1}^{n} x_i^2 \le 1 \},
\end{equation*}
where
\begin{equation}\label{setofcurve-1}
\Gamma_n^H: =\{ \gamma\in C(D_n, H)\ |\ \gamma(-x)=-\gamma(x),\ x\in D_n\ ; \ I(\gamma(x))<0,\ x\in \partial{D_n} \} ,
\end{equation}
and $H=E_r$ is the subspace of $E$, which consists of radial functions. The reason we write $c_n^H$ here is that in the future we will define several min-max values associated to different functional spaces $H$, and $H=E_r$ is only one of three choices.
By lemma~\ref{nonempty-1}, $\Gamma_n^H$ is not empty. We view $D_n$ as subset of $D_{n+1}$, so $\Gamma_{n+1}^H \subset \Gamma_n^H $ and ${c}_n^H$ is an increasing sequence.
By standard deformation theorem (see \cite{ST} or \cite{MW}), for each minimax value ${c}_n^H$, there is a corresponding Palais-Smale sequence $\{ u_i \}$ such that, as $i \to \infty$,
\begin{equation*}
I(u_i) \to {c}_n^H \ \ \ \textrm{and} \ \ \ I'(u_i) \to 0 \ \ \textrm{in $H^*$ strongly },
\end{equation*}
where $H^*$ is the dual space of $H$. Although we work on $E_r$, the compactness of the Palais-Smale sequence is not clear. Actually without Ambrosetti-Rabinowitz condition \cite{AR},
even the boundedness of this sequence is not obvious.
In order to get the boundedness and compactness of the Palais-Smale sequence, we have to choose a special Palais-Smale sequence with additional property.
Since any solution $u$ of \eqref{main} in $E$ satisfies Pohozaev identity $P(u)=0$ with
\begin{equation}\label{Pohozaev-2}
P(u):=(N-2)\int_{\RR^N} |\nabla{u}|^2dx-2N \int_{\RR^N} G(u)dx,
\end{equation}
we can construct a suitable extended functional $I_e(t,v)$ (see \eqref{extended-functional}),
on the augmented space $\RR \times H$ and manage to extract a palais-Smale sequence related to Pohozaev identity (see remark~\eqref{Palais-Smale-2}), as $i \to \infty$,
\begin{equation*}
I(u_i) \to c_n^H,   \ I'(u_i) \to 0 \ \ \textrm{in $H^*$ strongly }  \textrm{and} \  \ P(u_i) \to 0 .
\end{equation*}
This idea was first used in Jeanjean \cite{JEA} for a nonlinear eigenvalue problem and was fully developed in \cite{HIT}. By radial lemma~\ref{radial-lemma}, Egorov's theorem
and theorem~\ref{compactness-1} we obtain the compactness of this special Palais-Smale sequence. Overall we have the following existence results concerning radial solutions.
\begin{theorem}\label{existence}
Under (g1), (g2) and (g3), \eqref{main} has a sequence of radial solutions $\{\bar{u}_n\}$ and $I(\bar{u}_n)=c_n^{E_r}$. Moreover $\bar{u}_1$ is nonnegative least energy solution and $c_n^{E_r} \to \infty$ as $n \to \infty$.
\end{theorem}
When $g$ takes special form, Gidas, Ni and Nirenberg \cite{GNN1, GNN2} shows that positive solutions with certain decaying rate at infinity are always radially symmetric up to translation  by moving plane method.
One typical example for theorem~\ref{existence} is $g(s)=s^{p}-s^{q}$ for $ s \ge 0$ with $0 <q<p< 2^{*}-1$ for $N\ge 2 $. In particular when $0<q<1$, all solutions of \eqref{main} have compact supports, Cort\'azar, Elgueta and Felmer \cite{CEF1} obtain the existence and uniqueness of this nonnegative radial least energy solution up to translation.

To state the results about the existence of nonradial solutions, we need to introduce more notations. We follow the idea and notations in \cite{BW}. Assume that $N\ge 4$ and $2\le M\le N/2$.
Take $\tau \in \mathcal{O}(N)$ such that $ \tau(x^2,x^1,x^3)=(x^1,x^2,x^3) $ for $x^1\in \RR^M$, $x^2\in \RR^M$ and $x^3\in \RR^{N-2M}$, where $x=(x^1,x^2,x^3)\in \RR^N=\RR^M \times \RR^M \times \RR^{N-2M}$.
Define
\begin{equation*}
X_{\tau} :=\{ \ v\in E\ |\ v(\tau x)=-v(x) \ \textrm{ for all $x\in \RR^N$ } \}.
\end{equation*}
It is easy to see that $X_{\tau}$ contains no nontrivial radial functions and any nontrivial function in $X_{\tau}$ changes sign.
Let $E_{1}$ denote the subspace of invariant functions with respect to $\mathcal{O}_1=\mathcal{O}(M)\times \mathcal{O}(M) \times \textrm{id}\subset \mathcal{O}(N)$, which acts isometrically on $E$.
Let $E_{2}$ denote the subspace of invariant functions with respect to $\mathcal{O}_2=\mathcal{O}(M)\times \mathcal{O}(M) \times \mathcal{O}(N-2M) \subset \mathcal{O}(N)$, which acts isometrically on $E$. When $N=2M$, we agree that the components corresponding to $N-2M$ do not exist. It is clear that $E_2$ is a subspace of $E_1$ and they coincide when $N=2M$.

Let $H$ be $X_{\tau}\cap E_1$ or $X_{\tau}\cap E_2$ in \eqref{setofcurve-1}, repeating the technique of adding one dimension of space, we also can obtain a Palais-Smale sequence related to Pohozaev identity and manage to prove the boundedness of this sequence. So we have the following results about nonradial solutions.
\begin{theorem}\label{nonradial-existence}
Assume (g1), (g2) and (g3).
\begin{enumerate}
\item [(i)] If $N\ge 5$ and $N-2M\ge 1$, \eqref{main} has one nonradial solution $\tilde{u}\in X_{\tau}\cap E_1$ with $I(\tilde{u})>2c_1^{E_r}$, which minimizes the energy functional $I$ over all nontrivial solutions belonging to $X_{\tau}\cap E_1$.
\item [(ii)] If $N=4$ or $N\ge 6$, and $N-2M\neq 1$, \eqref{main} has a sequence of nonradial solutions $\{\hat{u}_n\}\subset X_{\tau}\cap E_2$ and $I(\hat{u}_n) \to \infty $ as $n \to \infty$. Moreover $I(\hat{u}_1)>2c_1^{E_r}$ and $\hat{u}_1$ minimizes the energy functional $I$ over all nontrivial solutions belonging to $X_{\tau}\cap E_2$.
\end{enumerate}
\end{theorem}
For $0<q<1$, by comparison principle from \cite{SQ5}, any solution of \eqref{main} in $E$ is compactly supported. From \cite{BJM}, any least energy solution of \eqref{main} is radially symmetric about some point and has constant sign. Then we can build high energy solution by putting several least energy solutions together with enough distances from the supports of each other. It is clear this high energy solution is not in $E_r$. If this solution is in $E_1$ or $E_2$, by symmetry and \cite{BJM}, each connected component of the support of this high energy solution is a ball with center $(0,0, x^3) \in \RR^N$ for some $x^3 \in \RR^{N-2M}$. Then this high energy solution has $\mathcal{O}(2M)$ symmetry, and it follows that this high energy solution is not in $X_{\tau}$. So the solutions produced by theorem~\ref{existence} and theorem~\ref{nonradial-existence} are nontrivial.

This paper is organized as follows. Preliminary analysis and Theorem~\ref{TRU-MOS} are done in section 2. Theorem~\ref{existence} and theorem~\ref{nonradial-existence} are established section 3, section 4 and section 5. In section 3 we prove part of theorem~\ref{existence} and theorem~\ref{nonradial-existence}, and show that $\{c_n^H\}$ is a sequence of critical values for $I(v)$. In section 4 we analyze the energy level $c_1^H$ and construct a special path to show that $c_1^H$ is the least energy level for solutions in $H$. Moreover $c_1^{E_r}$ is the least energy level for solutions in $E$ and $c_1^{E_r}= \inf \{ \ I(v) \ | \ P(v)=0, \ v \in E-\{0\} \ \}$, where $P$ is from \eqref{Pohozaev-2}. In order to show $\lim_{n \to \infty} c_n^H=\infty$,
we modify the equation \eqref{main} and get a comparison equation with Ambrosetti-Rabinowitz condition to achieve the goal, this is the work of section 5.
For simplicity, we write $(\int_{\RR^{2}} |v|^{p} dx)^{\frac{1}{p}}$ as $ ||v||_{p}$, and also write $e^t$ as $\exp(t)$. Noticing that $2^{*}=\infty $ if $N=2$ and $2^*=(2N)/(N-2)$ if $N\ge 3$, we also mean $2^{*}-1=\infty $ if $N=2$. $H$ may be one of these three functional spaces: $E_r$, $E_1\cap X_{\tau}$ and $E_2 \cap X_{\tau}$ and we will specify the space if necessary. From time to time, we will use  $C$, $C_1$, $C_2$, etc to denote some generic constants.

\section{Preliminaries}
In this section we will introduce some basic inequalities and estimates.
\begin{lemma}\label{inequality-1}
Let $l\ge 1$, $t\ge 0$, $s\ge 0$, $k\ge 0$ and $\alpha > 0$.
\begin{enumerate}
\item [(i)] $ (e^t-1)^l \le e^{tl} -1$.
\item [(ii)] Set $H_{\alpha}(t):=|t|^k(e^{\alpha t^2}-1)$.  Then $$H_{\alpha}(s+t)\le 2^{k-1}(H_{4\alpha}(s) +H_{4\alpha}(t))\ \ \textrm{and}\ \ H_{\alpha}(s)\le 2^{k-1}(H_{4\alpha}(s-t) +H_{4\alpha}(t))$$
\end{enumerate}
\end{lemma}
\begin{proof}
(i) is trivial since $l\ge 1$ and $ [(e^t-1)/e^t]^l +[1/e^t]^l \le 1$. Now for (ii), we see that by Taylor expansion the derivative of $H_{\alpha}'(t)$ is increasing, so $H_{\alpha}(t)$ is a convex function. Therefore by convexity
\begin{multline*}
H_{\alpha}(s+t)=|s+t|^k(\exp(\alpha (s+t)^2)-1)=2^{k}(s/2+t/2)^{k}(\exp(4\alpha(s/2+t/2)^2)-1)
 \\
=2^k H_{4\alpha}(s/2+t/2)=2^k((1/2)H_{4\alpha}(s)+(1/2)H_{4\alpha}(t))=2^{k-1}(H_{4\alpha}(s)+H_{4\alpha}(t)).
\end{multline*}
Therefore $ H_{\alpha}(s)\le  H_{\alpha}(|s-t|+t)\le  2^{k-1}(H_{4\alpha}(s-t) +H_{4\alpha}(t)) $
\end{proof}
Next lemma is due to Edmunds and Ilyin \cite{EI}, which is a sharp Gagliardo-Nirenberg type inequality on the plane $\RR^2$. This lemma is first part of theorem~\ref{TRU-MOS}. The proof given in \cite{EI} heavily replies on the embedding inequality $||v||_{(2p)/(2-p)} \le c(p) ||\nabla{v}||_p $ for any $v \in C_{0}^{\infty}(\RR^2) $. We will give a different proof to this inequality by exploiting the homogeneity of this inequality.
\begin{lemma}\label{inequality-2}
Let $v \in C_{0}^{\infty}(\RR^2) $ and $p \in [ q+1, \ \infty)$.  Then
\begin{equation}\label{EI-1}
||v||_{p} \le C(q,p) ||\nabla{v}||_{2}^{1-\frac{q+1}{p}} ||v||_{q+1}^{\frac{q+1}{p}} ,
\end{equation}
where the coefficient $C(q,p)=2(1+\Gamma(p/2 +1))^{1/p}$, $C(q,p)\le Q p^{1/2}$ and $Q$ is a constant depending only on $q$.
\end{lemma}
\begin{proof}
It is enough to prove that, for any radial function $u(x)=u(|x|)$ with $||\nabla{u}||_2 \le 1$, which is nonnegative, compactly supported and decreasing in $|x| \in [0, \infty)$, there holds
\begin{equation}\label{EI-2}
||u||_p^p \le (C(q,p))^p ||u||_{q+1}^{q+1}.
\end{equation}
Indeed rewrite \eqref{EI-1} into the form
\begin{equation*}
||v/||\nabla{v}||_{2}||_{p} \le C(q,p) ||v/||\nabla{v}||_{2}||_{q+1}^{\frac{q+1}{p}},
\end{equation*}
and let $u:=|v|/||\nabla{v}||_{2}$. By means of Schwarz symmetrization, we reach \eqref{EI-2}.

Now following Moser's idea in \cite{MOS}, see also \cite{AT}, set
\begin{equation*}
w(t):=\sqrt{4 \pi}u(e^{-t/2}),\ \ |x|^2=e^{-t}.
\end{equation*}
Then $w(t)$ is a nonnegative and increasing function defined on $(-\infty, \infty)$. Furthermore we have
\begin{equation*}
1\ge \int_{\RR^2}|\nabla{u}|^2dx =\int_{-\infty}^{\infty}|\frac{dw}{dt}|^2 dt,\ \ \int_{\RR^2}u^p dx=\frac{\pi}{(4 \pi)^{p/2}}\int_{-\infty}^{\infty} w^p e^{-t} dt.
\end{equation*}
Let $T_0:=\sup \{ t \ | \ w(t) \le 1 \} $. It is clear that $-\infty <T_0 \le \infty$. Observing that if $T_0=\infty$, we are done with $C(q,p)=1$. So we assume in the following that $T_0<\infty$. If $t \ge T_0$, by H\"older inequality, we have $w(t)= w(T_0)+\int_{T_0}^{\infty} (dw/dt ) dt\le 1+(t-T_0)^{1/2}$. Then we have, for $t\ge T_0$,
\begin{equation*}
0\le (w(t))^p \le (1+(t-T_0)^{1/2})^p \le 2^p(1+(t-T_0)^{p/2}).
\end{equation*}
Therefore by direct calculation, we get
\begin{eqnarray*}
 & & ||u||_p^p = \frac{\pi}{(4\pi)^{p/2}}\int_{-\infty}^{\infty}w^p(t)e^{-t} dt= \frac{\pi}{(4\pi)^{p/2}}[\int_{-\infty}^{T_0}w^p(t)e^{-t} dt+\int_{T_0}^{\infty}w^p(t)e^{-t} dt] \\
  &\le &\frac{\pi}{(4\pi)^{p/2}}[\int_{-\infty}^{T_0}w^{q+1}(t)e^{-t} dt+2^p \int_{T_0}^{\infty}(1+(t-T_0)^{p/2})e^{-t} dt] \\
   &= &\frac{\pi}{(4\pi)^{p/2}}[\int_{-\infty}^{T_0}w^{q+1}(t)e^{-t} dt+2^p (1+\Gamma(1+p/2)) \int_{T_0}^{\infty} e^{-t} dt] \\
   &\le &\frac{\pi}{(4\pi)^{p/2}}[\int_{-\infty}^{T_0}w^{q+1}(t)e^{-t} dt+2^p (1+\Gamma(1+p/2)) \int_{T_0}^{\infty}w^{q+1} e^{-t} dt] \\
   &\le &\frac{\pi}{(4\pi)^{p/2}}2^p (1+\Gamma(1+p/2))\int_{-\infty}^{\infty}w^{q+1}(t)e^{-t} dt
\end{eqnarray*}
Finally we reach
\begin{equation*}
  ||u||_p^p \le  \frac{ 2^p (1+\Gamma(1+p/2))}{(4\pi)^{(p-q-1)/2}}[\int_{-\infty}^{\infty}\frac{\pi}{(4\pi)^{(q+1)/2}}w^{q+1}(t)e^{-t} dt]\le [C(q,p)]^p||u||_{q+1}^{q+1},
\end{equation*}
where the coefficient $C(q,p)=2(1+\Gamma(p/2 +1))^{1/p}$. With the help of Stirling's asymptotic formula for gamma function, we have $C(q,p)\le Q p^{1/2}$ and $Q$ is a constant depending only on $q$.
\end{proof}
By lemma~\ref{inequality-2}, for $N=2$, Banach space $E$ can be embedded into any $L^{p}(\RR^2)$ for any $p\ge q+1$.
Combining lemma~\ref{inequality-2} with Taylor expansion, we obtain the following theorem, which includes the second part of theorem~\ref{TRU-MOS}.
\begin{theorem}\label{TRU-MOS-F}
Let $N=2$. There exists a small constant $\rho_1>0$ depending only on $q$ such that, for any nonnegative number $k$ with $k+1\ge q$,
 $$ \int_{\RR^2}|v|^{k}(e^{v^2} -1)dx \le C ||\nabla{v}||_{2}^{k+1-q}||v||_{q+1}^{q+1}  \ \ \ \   \textrm{for any $v\in E$ with  $||\nabla{v}||_2 \le \rho_1 $ ,}$$
where $C$ depends on $k$ and $q$. In particular, if $k+1\ge q$ and $||v||_E \to 0$, then $\int_{\RR^2}|v|^k(e^{v^2} -1)dx \to 0$.
\end{theorem}
\begin{proof}
By Taylor expansion formula, we have that $e^{t^2}-1=\sum_{p=1}^{\infty}t^{2p}(p!)^{-1}$. So by lemma~\ref{inequality-2}, we have
\begin{eqnarray*}
 & &\int_{\RR^2}|v|^k(e^{v^2}-1)dx = \int_{\RR^2}\sum_{p=1}^{\infty}v^{2p+k}(p!)^{-1}dx=\sum_{p=1}^{\infty}\int_{\RR^2}v^{2p+k}(p!)^{-1}dx \\
 &=&\sum_{p=1}^{\infty}(p!)^{-1}||v||_{2p+k}^{2p+k}=\sum_{p=1}^{\infty}(p!)^{-1}Q^{2p+k}(2p+k)^{(2p+k)/2} ||\nabla||_2^{2p+k-q-1}||v||_{q+1}^{q+1} \\
 &=&(||\nabla{v}||_{2}^{k+1-q}||v||_{q+1}^{q+1})\sum_{p=1}^{\infty}(p!)^{-1}Q^{2p+k}(2p+k)^{(2p+k)/2} ||\nabla{v}||_2^{2p-2}.
\end{eqnarray*}
For the series $\sum_{p=1}^{\infty}(p!)^{-1}Q^{2p+k}(2p+k)^{(2p+k)/2} ||\nabla{v}||_2^{2p-2}$, the ratio of latter term over previous term is
$$Q^2 \frac{2p+2+k}{p+1}(1+\frac{2}{2p+k})^{\frac{2p+k}{2}}||\nabla{v}||_2^2.$$
As $p \to \infty$, the limit is $2eQ^2||\nabla{v}||_2^2$. Therefore taking $\rho_1=\frac{1}{2} Q^{-1}(\sqrt{e})^{-1}$, we get the desired inequality. If $||v||_E \to 0$, we have $||\nabla{v}||_2+||v||_{q+1} \to 0$.
By the inequality we just developed, we get $\int_{\RR^2}|v|^k(e^{v^2} -1)dx \to 0$.
\end{proof}
Recall the energy functional associated with \eqref{main} is $I(v)= \frac{1}{2} \int_{\RR^{N}} |\nabla{v}|^{2}dx- \int_{\RR^{N}}G(v)dx$, for any $ v \in E$, and the Banach space
$E$ is the span of $ C_0^{\infty}(\RR^{N})$ under the norm $ ||v||_{E} = ||\nabla v||_2 + ||v||_{q+1} $. Next lemma says that $I(v)$ is well defined.
\begin{lemma}\label{makesense}
Assume $(g1)$ and $(g2)$. Then, for any $v\in E$, $\int_{\RR^{N}}G(v)dx$ is finite. Moreover there exists a small constant $\rho_2>0$ such that, for any $v\in E$ with  $||\nabla{v}||_2 \le \rho_2 $,
\begin{eqnarray*}
\int_{\RR^2}G(v)dx \le - \frac{m}{q+1}||v||_{q+1}^{q+1}  \ \ \textrm{if  $N=2 $,} \ \ \int_{\RR^N}G(v)dx \le - \frac{m}{q+1}||v||_{q+1}^{q+1}+ \frac{||\nabla{v}||_2^2}{4}  \ \   \textrm{if  $N\ge 3 $,}
\end{eqnarray*}
where $\rho_2$ depends on $q$ and $m$.
\end{lemma}
\begin{proof}
First we deal with the case $N\ge 3 $. By  $(g1)$ and $(g2)$, we have that $|G(s)|\le C(|s|^{q+1}+|s|^{2^*})$. So it is clear that $\int_{\RR^{N}}G(v)dx$ is finite. Using the growth assumptions again, with the help of Sobolev embedding inequality, we find that
$$\int_{\RR^N} G(v) dx \le \int_{\RR^N} -\frac{m}{q+1}|v|^{q+1}+ C|v|^{2^*}dx \le -\frac{m}{q+1} ||v||_{q+1}^{q+1}+C||\nabla{v}||_2^{2^*} .$$
Choosing $\rho_2=\min\{ 1, (1/(4C))^{(N-2)/4}\}$, we can get the desired estimate.

Next we deal with the case $N=2 $. By  $(g1)$ and $(g2)$,  for any $\alpha >0$, we have
\begin{equation}\label{inequality-4}
|g(s)|\le C(|s|^q+|s|^q (e^{\alpha s^2}-1))\  \ \textrm{and} \ \   |G(v)| \le C( |v|^{q+1}+|v|^{q+1} (e^{\alpha v^2}-1)).
\end{equation}
For any nonzero $v$, choose $\alpha$ such that $\sqrt{\alpha}||\nabla{v}||_2 =\rho_1$. From theorem~\ref{TRU-MOS-F} and the choice of $\alpha$, we see that $\int_{\RR^2}|G(v)|dx $ is finite.
By $(g1)$ and $(g2)$ again,  $g(s)\le -\frac{3m}{2}s^q+Cs^q(e^{s^2}-1)$ for $s \ge 0$, we find that $G(v)\le
-\frac{3m}{2(q+1)} |v|^{q+1}+C|v|^{q+1} (e^{v^2}-1)$. From theorem~\ref{TRU-MOS-F}, we get $ \int_{\RR^2}G(v)dx \le (-\frac{3m}{2(q+1)}+C||\nabla{v}||_2^{2})||v||_{q+1}^{q+1}\le -\frac{m}{q+1} ||v||_{q+1}^{q+1}$
by choosing $\rho_2=\min\{\rho_1, \sqrt{\frac{m}{2C(q+1)} }  \}$, where $\rho_1$ is from theorem~\ref{TRU-MOS-F}.
\end{proof}
Next we show that $I(v)$ is continuously differentiable from $ E $ to $\RR$. Recall from lemma~\ref{inequality-1}, $H_{\alpha }(t)=|t|^k(e^{ \alpha t^2}-1)$.
\begin{theorem}\label{c1functional}
Assume $(g1)$ and $(g2)$.

(i) For $N=2$,  any $k\ge q-1$ and any $v \in E$,  $\int_{\RR^2}H_1(v)dx < \infty$ and both functionals $\int_{\RR^2}H_1(v)dx$ and $\int_{\RR^2}G(v)dx $ are continuous from $E$ to $\RR$.

(ii) For $N\ge 2$, $I(v):E \to \RR$ is of $C^{1}$.
\end{theorem}
\begin{proof}
For $N\ge 3 $, following the arguments in \cite{BL1},  it is quite standard to show that $I(v):E \to \RR$ is of $C^{1}$. Therefore we prove this theorem only for the case $N=2$.
For any $v\in E$, by density we can find $w \in C_{0}^{\infty}(\RR^2)$ such that $2||\nabla{(v-w)}||_2 \le \rho_1$, where $\rho_1$ is from theorem~\ref{TRU-MOS-F}.
Therefore $\int_{\RR^2}H_4(v-w)dx < \infty$ and $\int_{\RR^2}H_4(w)dx < \infty$.
By lemma~\ref{inequality-1} we have
\begin{equation}\label{inequality-3}
H_{1}(v)=H_{1}(|v|)\le 2^{k-1}(H_4(|v-w|)+H_4(|w|),
\end{equation}
which simply means that $\int_{\RR^2}H_1(v)dx < \infty$ is finite. Now for any fixed $w \in E$, let $v_i$ goes to $w$ strongly in $E$ as $i \to \infty$, taking subsequence if necessary,
we may assume that $v_i $ converges to $w$ pointwise a.e. in $\RR^2$. From theorem~\ref{TRU-MOS-F}, $H_4(|v_i-w|)$ converges to $0$ in $L^1(\RR^2)$. By \eqref{inequality-3}, $H_{1}(v_i)\le 2^{k-1}(H_4(|v_i-w|)+H_4(|w|))$.
By Lebesgue convergence theorem, as $i$ goes to $\infty$, $\int_{\RR^2} H_{1}(v_i)dx$ converges to $\int_{\RR^2} H_{1}(w)dx$, which implies the continuity. Using \eqref{inequality-4} with $k=q+1$ and Lebesgue convergence theorem, we have the continuity of $J $ from $E$ to $\RR$, where $J(v):= \int_{\RR^2}G(v)dx $.

Now in order to show that $I(v)$ is continuously differentiable, from the structure of $I(v)$, we only need to show that $J(v)$ is
continuously differentiable.  Let $u, \ v \in E$, define $K(v):= \int_{\RR^2}g(u)vdx $. From $(g1)$, $(g2)$, \eqref{inequality-4} and lemma~\ref{inequality-1}, we have
\begin{equation}\label{inequality-5}
|g(u)|^{\frac{q+1}{q}}\le C(|u|^{q+1}+|u|^{q+1}(e^{\frac{q+1}{q}u^2}-1))\ \ \textrm{and} \ \ \int_{\RR^2}|g(u)|^{\frac{q+1}{q}}dx<\infty.
\end{equation}
By H\"{o}lder inequality, we see that $K(v)$ is a bounded linear functional on $E$.
Since $g$ is continuous, for $|t| \le 1$, by \eqref{inequality-4}, \eqref{inequality-5} and mean value theorem we see that
$$|\frac{G(u+tv)-G(u))}{t}| \le C \{(|u|+|v|)^q+(|u|+|v|)^q(e^{(|u|+|v|)^2-1})\}|v| \in L^{1}(\RR^2).$$
From Lebesgue convergence theorem, we see that the functional $J$ has the G\^{a}teaux differential $K$ at $u$.
Let $u_i \in E$ such that $u_i \to u$ strongly in $E$ as $i \to \infty$, from \eqref{inequality-5} and Lebesgue convergence theorem, we find that $\int_{\RR^2}|g(u_n)-g(u)|^{\frac{q+1}{q}}dx \to 0$ as $i \to \infty$.
Finally we have $|\int_{\RR^2}(g(u_i)-g(u))v dx| \le ||g(u_i)-g(u)||_{\frac{q+1}{q}} ||v||_{q+1}\le ||g(u_i)-g(u)||_{\frac{q+1}{q}} ||v||_{E}$. Therefore we conclude that
$$ \sup \{ |\int_{\RR^2} (g(u_i)-g(u))vdx | \ \ | \ v \in E \ , \ ||v||_E \le 1     \} \to 0 \textrm{ as $i \to \infty$}, $$
which implies $J(v)$ is $C^1$ from $E$ to $\RR$. So is $I(v)$ .
\end{proof}
Recall that $E_r$ is a subspace of $E$, which consists of radial functions. In order to develop some compactness results, we give the following radial lemma
concerning the uniform decay rate of radial functions in $E_r$ at infinity. When $q=1$,  this lemma is due to Strauss \cite{STR}, see also \cite{BL1}.
\begin{lemma}\label{radial-lemma}
There exists a constant $C$, depending on $q$, such that for all $v \in E_r$ there holds
$$ |v(x)|\le (\frac{q+3}{2\omega_{N-1}})^{\frac{2}{q+3}} \frac{||v||_E}{|x|^{\frac{2}{q+3}(N-1)}}, \ \ \forall\  |x|\ge 1 ,$$
where $\omega_{N-1}$ is the area of the unit sphere in $\RR^N$.
\end{lemma}
\begin{proof}
For simplicity we write $|x|$ as $r$ and denote $v(x)$ as $v(r)$. Let $\beta:=\frac{q+3}{4}$, taking any $v(r) \in C_{0}^{\infty}(\RR^N)\cap E_r$, we may assume $v\ge 0$,
otherwise we use $|v(r)|$ to replace $v$ since $|v|$ is Lipschitz continuous. Therefore we have
\begin{equation*}
\frac{d(r^{N-1}(v(r))^{2\beta})}{dr} \ge 2\beta r^{N-1} (v(r))^{2\beta-1}\frac{d(v(r))}{dr}.
\end{equation*}
For $r \ge 1$, integrating over $ [r,\ \infty) $, by H\"{o}lder inequality we have
$$(v(r))^{2\beta}r^{N-1} \le - \frac{2\beta}{\omega_{N-1}}\int_{r}^{\infty}\frac{d(v(s))}{ds}{\omega_{N-1} s^{N-1}} (v(s))^{2\beta-1} ds\le \frac{2\beta}{\omega_{N-1}} ||v||_{4\beta-2}^{2\beta-1}||\nabla{v}||_2.$$
Noticing that $\beta=\frac{q+3}{4}$, we have $(v(r))^{2\beta}r^{N-1} \le \frac{2\beta}{\omega_{N-1}} ||v||_{q+1}^{2\beta-1}||\nabla{v}||_2\le \frac{2\beta}{\omega_{N-1}} ||v||_{E}^{2\beta}.$
\end{proof}
For the function $g(s)$, satisfying (g1) and (g2), we define $g_{1}(s):=\max\{g(s)+ms^q, 0\}$ if $s\ge 0$ and
$g_{1}(s):=-g_{1}(-s)$ if $s\le 0$, and define $G_1(s):=\int_{0}^{s}g_1(t)dt$. In viewing of (g1), there exists small $\delta>0$ such that $g_1(s)=G_1(s)=0$ if $|s| \le \delta$.
Next we give some compactness results, part $(i)$ is the classical compactness results for $N\ge 2$ and part $(ii)$ will be used in the following sections.
\begin{theorem}\label{compactness-1}
Assume (g1), (g2), and $v_i \rightharpoonup w$ weakly in $E_r$ or $ E_2$, and pointwise a.e. in $\RR^N$ as $i \to \infty$. $N=4$ or $N\ge 6$, and $N-2M\neq 1$ are assumed for $E_2$.
\begin{enumerate}
\item [(i)] If $q+1<p<2^*$, then $\lim_{i \to \infty} ||v_i-w||_p =0$, that is , the embedding $E_r\hookrightarrow L^{p}(\RR^N)$ and $ E_2 \hookrightarrow L^{p}(\RR^N)$ are compact.
\item [(ii)]
There hold $\int_{\RR^N}G_1(v_i)dx \to \int_{\RR^N}G_1(w)dx $ and $\int_{\RR^N}g_1(v_i)v_idx \to \int_{\RR^N}g_1(w)wdx$ as $i $ goes to $\infty$.
\end{enumerate}
\end{theorem}
\begin{proof}
By weak convergence we have that $||v_i||_{E} \le M^*$ for some $M^*>0$.
(i) We first prove that $E_r\hookrightarrow L^{p}(\RR^N)$ is compact. From lemma~\ref{inequality-2} and Sobolev embedding inequality for $N\ge 3$, $||v_i||_p \le CM^*$. By linearity we only need to show that this is true for $w=0$. Choose  $\theta >1$ such that $\theta p <2^*$. Let $\theta' = \theta/(\theta-1) $. By lemma~\ref{radial-lemma}, for any $\epsilon>0$, there exists $R_0>1$ such that $ |v_i(x)|^{p-q-1}\le \epsilon $ when $|x|\ge R_0$. Hence
 $$\int_{\RR^N-B(0,R_0)} |v_i|^{p} dx \le \epsilon \int_{\RR^N-B(0,R_0)} |v_i|^{q+1} dx \le \epsilon (M^*)^{q+1}. $$
Since $v_i $ converges to $0$  pointwise a.e. in $B(0,R_0)$, by Egorov's theorem, there exists a measurable subset $U_\epsilon$ with Lebesgue measure $|U_\epsilon| \le \epsilon^{\theta'}$
such that $v_i(x) $ converges to $0$ uniformly in $(B(0,R_0)-U_\epsilon)$. Then there exists $i_\epsilon>0$ such that, if $i \ge i_\epsilon$, we have $$\int_{B(0,R_0)-U_\epsilon} |v_i|^{p} dx \le \epsilon.$$
Finally by H\"{o}lder inequality, we have $$\int_{U_\epsilon}|v|_i^{p} dx \le ||v_i||_{\theta p}^{p}|U_\epsilon|^{1/\theta'} \le C (M^*)^{p}\epsilon.$$
So for $i \ge i_\epsilon$, we have $\int_{\RR^N}|v|_i^{p} dx \le (1+(M^*)^{q+1}+C(M^{*})^p)\epsilon$, which implies $\lim_{i \to \infty} ||v_i||_p =0$.

Next we show that $E_2 \hookrightarrow L^{p}(\RR^N)$ is compact. Since in this case $\RR^N$ is compatible with $\mathcal{O}_2$, see \cite{MW}, there exists some $r_0>0$ such that
\begin{equation}\label{ball-1}
\lim_{i \to \infty} \sup_{y\in \RR^N}\int_{B(y,r_0)}|v_i|^{q+1} dx =0.
\end{equation}
By lemma I.1. from \cite{LION}, we have $\lim_{i \to \infty} ||v_i||_p =0$.

(ii) Let $g_N(s):=|s|^{q+1}(e^{\alpha s^2}-1) $ with $\sqrt{2\alpha}M^*=\rho_1$ for $N=2$ and $g_N(s):=|s|^{2^*}$ for $N \ge 3$. Take some $p$ with $q+1<p<2^*$. For any $\epsilon>0$, by (g2) and definition of $g_1$, we have $G_{1}(v_i)+g_1(v_i)v_i \le C |v_i|^p+ \epsilon g_N(v_i) $. By Fatou's lemma, we have
\begin{equation*}
\liminf_{i \to \infty}\int_{\RR^N} C |v_i|^p+ \epsilon g_N(v_i)-G_{1}(v_i) dx \ge \int_{\RR^N} C |w|^p+ \epsilon g_N(w)-G_{1}(w) dx.
\end{equation*}
But from (i), $\lim_{i \to \infty}\int_{\RR^N} |v_i|^p dx =\int_{\RR^N} |w|^p dx $, and $ \int_{\RR^N} g_N(v_i)dx $ is uniformly bounded,
so we have $\limsup_{i \to \infty} \int_{\RR^N}G_1(v_i)dx \le \int_{\RR^N}G_1(w)dx $. From Fatou's lemma again,
$$\liminf_{i \to \infty} \int_{\RR^N}G_1(v_i)dx \ge \int_{\RR^N}G_1(w)dx , $$
which means $\lim_{i \to \infty} \int_{\RR^N}G_1(v_i)dx = \int_{\RR^N}G_1(w)dx $. Likewise we can get $\lim_{i \to \infty} \int_{\RR^N}g_1(v_i)v_idx$ $= \int_{\RR^N}g_1(w)wdx$ .
\end{proof}
\begin{remark}
We can show that $E_r \hookrightarrow L^{p}(\RR^N)$ is compact without using lemma~\ref{radial-lemma}. Indeed $\RR^N$ is also compatible with $\mathcal{O}(N)$, so we have \eqref{ball-1}. By lemma I.1. from \cite{LION}, we have the compactness. Therefore lemma~\ref{radial-lemma} is not the only path to get compactness.
\end{remark}
\section{Mountain pass structure and minimax values}
To the end of this paper, we always assume that (g1),(g2) and (g3) hold. From previous section we know that $I(v)$  is continuously differentiable from $E$ to $\RR$. We say $u$ is a critical point of $I(v)$ if $I'(u)=0$, that is
\begin{equation*}
\int_{\RR^N}\nabla{u} \nabla{v} dx-\int_{\RR^N}g(u) v dx =0, \ \ \forall \ v \in E.
\end{equation*}
This exactly says that $u$ is a weak solution of \eqref{main}. Therefore finding solutions of \eqref{main} becomes searching for critical points of $I(v)$. We first present a result from \cite{BL2},
concerning the existence of an odd continuous mapping from $\partial{D_n}$ to $E_r(R)$, where $E_r(R)$ consists of radial functions in $E_r$ with supports inside the ball $B(0,R)$, $D_n=\{ x \in \RR^n \ |\ |x| \le 1 \}$ and
$\partial{D_n}=\{ x  \in D_n \ |\ |x| = 1 \}$.
\begin{lemma}\label{odd-mapping}
For all $n \ge 1$, there exist two positive constants $C_1(n)$, $C_2(n)$ depending on $n$, and an odd continuous mapping $\pi_n:\ \partial{D_n}\rightarrow E_r(R)$ such that $0 \notin \pi_{n}(\partial{D_n})$ and
$$ \int_{\RR^N} G(v) dx \ge C_1(n)R^N-C_2(n)R^{N-1} ,  \ \ \forall R> n. $$
Moreover function $v$ is Lipschitz continuous and $||v||_\infty  =s_0$ for any $v\in \pi_{n}(\partial{D_n}) $, where $s_0 $ is from (g3).
\end{lemma}
Next we show that $I(v)$ has mountain pass structure and symmetric mountain pass structure around $v=0$ and $v=0$ is a local minimizer.
\begin{lemma}\label{mountain-pass-structure-1}
\begin{enumerate}
\item [(i)] There exist $a_3>0$ and $\rho_3>0$ such that $$I(v)\ge 0, \ \forall \ ||v||_E \le \rho_3 ; \ \ I(v)\ge a_3 , \  \forall \ ||v||_E = \rho_3.  $$
\item [(ii)] For all $n \ge 1$, there exist two odd continuous mapping $\bar{\gamma}_n:\ \partial{D_n}\rightarrow E_r$ and $\hat{\gamma}_n:\ \partial{D_n}\rightarrow E_2\cap X_\tau $ such that $I(\bar{\gamma}_n(x))<0$ and $I(\hat{\gamma}_n(x))<0$ for all $x\in \partial{D_n}$.
\end{enumerate}
\end{lemma}
\begin{proof}
(i)\ Let $p=\max(2, q+1)$ and $\delta=\min (\frac{1}{4},\frac{m}{q+1})$. By lemma~\ref{makesense}, for any $v\in E$ with  $||v||_E \le \min(\rho_2,1) $,
\begin{equation*}
I(v)\ge \frac{1}{4}||\nabla{v}||_{2}^2 + \frac{m}{q+1}||v||_{q+1}^{q+1} \ge \delta(||\nabla{v}||_2^p+||v||_{q+1}^{p})\ge \frac{\delta}{2^{p-1}} ||v||_E^{p},
\end{equation*}
which implies $\rho_3=\min(\rho_2,1)$ and $a_3=\frac{\delta}{2^{p-1}} \rho_3^{p}$.

(ii) The construction of $\bar{\gamma}_n$ is straightforward. For fixed $n$, from lemma~\ref{odd-mapping}, take a large $R>n$ such that the odd continuous mapping $\pi_n:\ \partial{D_n}\rightarrow E_r(R)$ satisfies
\begin{equation*}
\int_{\RR^N} G(v) dx \ge C_1(n)R^N-C_2(n)R^{N-1} \ge 1.
\end{equation*}
By compactness we assume $||\nabla{\pi_n(x)} ||_2 \le C_3(n) $ for some $C_3(n)>0$ and any $x\in \partial{D_n}$. Let $l_n:=C_3(n)$ and define $\bar{\gamma}_n(x):=\pi_n(x)(y/l_n)$ for any $x\in \partial{D_n}$ and $y \in \RR^N$.
By calculation $$I(\bar{\gamma}_n(x))=\frac{1}{2}l_n^{N-2}||\nabla{\pi_n(x)}||_{2}^2-l_n^N\int_{\RR^N}G(\pi_n(x))dz \le \frac{1}{2}C_3^N(n)-C_3^N(n)=-\frac{1}{2}C_3^N(n)<0.$$

Next we construct $\hat{\gamma}_n$. Let $A:=\max_{|s|\le s_0}|G(s)|$. Take a smooth function $\phi(s)$ such that $\phi(s)=1 $ for $s \ge 1$, $\phi(s)=-1 $ for $s \le -1$ and $|\phi(s)| \le 1 $ for $-1\le s \le 1$. For any $u\in \pi_{n}(\partial{D_n})$, where $\pi_{n}$ is from the construction of $\bar{\gamma}_n$. Let us define
\begin{equation*}
\hat{u}(x^1, x^2, x^3):=u(x^1, x^2, x^3)\phi(|x^1|-|x^2|),\ \ \textrm{for}\  x^1\in \RR^M,\  x^2\in \RR^M \ \textrm{and}\  x^3\in \RR^{N-2M}.
\end{equation*}
It is clear $\hat{u} \in E_2\cap X_\tau $. We write $|x|=|(x^1,x^2,x^3)|$, $|x^1|$, $|x^2|$ and $|x^3|$ as $r$, $r_1$, $r_2$ and $r_3$ respectively. It is obvious that $r^2=r_1^2+r_2^2+r_3^2$. By direct calculation we have
\begin{eqnarray*}
& &\int_{\RR^N} G(\hat{u})dx=\int_{0}^{R} \int_{0}^{R} \int_{0}^{R+1} G(\hat{u}) r_1^{M-1} r_2^{M-1} r_3^{N-2M-1} dr_1 dr_2 dr_3 \\
&=&2 \int_{0}^{R} \int_{0}^{R} \int_{r_2}^{R+1} G(\hat{u}) r_1^{M-1} r_2^{M-1} r_3^{N-2M-1} dr_1 dr_2 dr_3 \\
&= & 2 \int_{0}^{R} \int_{0}^{R} \int_{r_2}^{r_2+1}( G(\hat{u})-G(u)) r_1^{M-1} r_2^{M-1} r_3^{N-2M-1} dr_1 dr_2 dr_3 \\
&+ & 2 \int_{0}^{R} \int_{0}^{R} \int_{r_2}^{R+1} G(u) r_1^{M-1} r_2^{M-1} r_3^{N-2M-1} dr_1 dr_2 dr_3.
\end{eqnarray*}
Therefore, choosing $R$ large enough, from lemma~\ref{odd-mapping} we have
\begin{eqnarray*}
& &\int_{\RR^N} G(\hat{u})dx \ge \int_{\RR^N} G(u)dx-\frac{4AR^{N-2M}}{M(N-2M)}\int_{0}^{R} [(r_2+1)^M-r_2^M] r_2^{M-1} dr_2 \\
&  \ge & \int_{\RR^N} G(u)dx-C(N)R^{N-1}\ge C_1(n)R^N-C_2(n)R^{N-1}-C(N)R^{N-1} \ge 1.
\end{eqnarray*}
Remember that when $N=2M$, the components corresponding to $N-2M$ do not exist.
Then taking $\hat{\pi}_{n}(y):=\phi(r_1-r_2) (\pi_{n}(y) )$ for all $y \in \partial{D_n}$, repeat the construction process for $\bar{\gamma}_n$ and we get the desired odd continuous mapping $\hat{\gamma}_n:\ \partial{D_n}\rightarrow E_2\cap X_\tau $.
\end{proof}
\begin{remark}\label{structure-1}
Since $E_2$ is a subspace of $E_1$, lemma~\ref{mountain-pass-structure-1} shows that, for all $n \ge 1$, the odd continuous mapping $\hat{\gamma}_n$ is also from $ \partial{D_n} $ to $ E_1 \cap X_\tau $ with $I(\hat{\gamma}_n(x))<0$ for all $x\in \partial{D_n}$.
\end{remark}
Now take $H$ to be any one of these three functional spaces: $E_r$, $E_1\cap X_{\tau}$ and $E_2 \cap X_{\tau}$. As stated in the introduction, we define a sequence of minimax values $c_n^H$:
\begin{equation}\label{minimax-1}
c_{n}^H:=\inf_{\gamma \in \Gamma_{n}^H} \sup_{x\in D_n} I(\gamma(x)),
\end{equation}
where $\Gamma_n^H:=\{ \gamma\in C(D_n, H)\ |\ \gamma(-x)=-\gamma(x),\ x\in D_n\ ; \ I(\gamma(x))<0,\ x\in \partial{D_n} \}. $
\begin{lemma}\label{nonempty-1}
For all $n \ge 1$, $\Gamma_n^H \neq \emptyset $ and $c_{n+1}^H \ge c_n^H \ge a_3$, where $a_3 $ is from lemma~\ref{mountain-pass-structure-1}.
\end{lemma}
\begin{proof}
Let $H=E_r$. From lemma~\ref{mountain-pass-structure-1}, set $\gamma_n(x):=|x|\bar{\gamma}_n(x/|x|) $ if $x\in D_{n}-\{0\}$ and $\gamma_n(0):=0$. It is clear that $\gamma_n$ is continuous and $\gamma_n\in \Gamma_n^{H}$. Let $x \in D_n$, then $(x,0)\in D_{n+1}$.
Through this identification, we view $D_n$ as a subset of $D_{n+1}$ and we have $\Gamma_{n+1}^H \subset \Gamma_n^H $, so we have $c_{n+1}^H \ge c_n^H  $.
By continuity we find that $\gamma(D_n) \cap \{ v\in E \ |\ ||v||_E =\rho_3\}\neq \emptyset$ for any $\gamma \in \Gamma_n^H$, this means that $c_n^H \ge a_3$. With lemma~\ref{mountain-pass-structure-1}, we can treat the other two cases almost identically.
\end{proof}
From standard deformation theorem (see \cite{ST}), for each minimax value $c_n^H$, we can select a Palais-Smale sequence $\{ u_i \}\subset H $ such that, as $i \to \infty$
\begin{equation*}
I(u_i) \to c_n^H \ \ \ \textrm{and} \ \ \ I'(u_i) \to 0 \ \ \textrm{in $H^*$ strongly },
\end{equation*}
where $H^*$ is the dual space of $H$. Unfortunately we are unable to prove the boundedness and compactness of this Palais-Smale sequence
due to lack of more structure information from $g$, such as Ambrosetti-Rabinowitz condition \cite{AR}. We need to find a special Palais-Smale sequence with more information.
Observing that any solution $u$ of \eqref{main} in $E$ satisfy Pohozaev identity
\begin{equation*}
P(u)=(N-2)\int_{\RR^N} |\nabla{u}|^2dx-2N \int_{\RR^N} G(u)dx=0,
\end{equation*}
we define the following extended functional $I_{e}:\ \RR \times E \rightarrow \RR$:
\begin{equation}\label{extended-functional}
I_e(t,v):=\frac{e^{(N-2)t}}{2}||\nabla{v}||_{2}^2-e^{Nt}\int_{\RR^N} G(v)dx,
\end{equation}
where $\RR \times E$ is equipped with standard product norm $||(t,v)||:=|t|+||v||_E$. With this norm, $I_e(t,v)\in C^{1}(\RR \times E, \RR)$ and the full derivative of $I_e$ is
$I_e'(t,v):=(( I_e)_t'(t,v), (I_e)_v'(t,v))$, where $( I_e)_t'(t,v)$ and $(I_e)_v'(t,v))$ are partial derivatives of $I_e(t,v)$ with respect to $t$ and $v$ respectively. There is an interesting relation between $I(v)$ and $I_e(t,v)$:
\begin{equation}\label{relation}
I_e(t,v)=I(v(\frac{x}{e^t})), \ \ \forall \ t\in \RR, \  \forall \ v\in E.
\end{equation}
Corresponding to \eqref{minimax-1}, we define another sequence of minimax values $b_n^H$:
\begin{equation}\label{minimax-2}
b_{n}^H:=\inf_{\lambda \in \Lambda_{n}^H} \sup_{x\in D_n} I_e(\lambda(x)),
\end{equation}
where $\Lambda_n^H:=\{ \lambda=(\tau(x), \gamma(x)) \in C(D_n, \RR \times H)\ |\ (\tau(-x),\gamma(-x))=(\tau(x), -\gamma(x)),\ x\in D_n\ ; \ I_e(\tau(x), \gamma(x))<0,\ x\in \partial{D_n} \}. $
\begin{lemma}\label{coincide-1}
For all $n \ge 1$, $\Lambda_n^H \neq \emptyset $ and $c_n^H =b_n^H$.
\end{lemma}
\begin{proof}
From lemma~\ref{nonempty-1}, we know that $\Gamma_n^H \neq \emptyset$. Taking any $\gamma \in \Gamma_n^H$, it is easy to see that $(0,\gamma_n)\subset \Lambda_n^H$, so $\Lambda_n^H \neq \emptyset $ and $c_n^H\ge b_n^H$.
On the other hand, taking any $\lambda(x)=(\tau(x), \widetilde{\gamma}(x))\in \Lambda_n^H$, let $\gamma(x)=\widetilde{\gamma}(x)(e^{-\tau(x)}y)$, it is clear that $I(\gamma(x))=I_e(\tau(x),\widetilde{\gamma}(x))$.
So we have $\gamma \in \Gamma_n^H$ and $c_n^H \le b_n^H$, which implies $c_n^H=b_n^H$.
\end{proof}
The advantage of using the extended functional $I_e(t,v)$ will be clear in the lemma below, which provides an additional condition related to Pohozaev identity. The distance between a point and a set is defined by
$ \textrm{dist}_{\RR \times H}((t_1, v_1), A)= \inf_{(t,v)\in A} (|t_1-t|+||v_1-v||_E)$ for any $A \subset \RR \times H$.
\begin{lemma}\label{coincide-2}
For each miniax value $c_n^H$, there exists a sequence $\{(t_i, u_i)\} \subset \RR \times H$ such that, as $i \to \infty$,
\begin{eqnarray}\label{Palais-Smale-1}
 & & t_i \to 0,  \  I_e(t_i,u_i) \to c_n^H,  \ (I_e)_v'(t_i,u_i) \to 0 \ \textrm{strongly in }\ H^{*}, \\ \nonumber
& &  P(u_i(\cdot/e^{t_i}))=\frac{(N-2)e^{(N-2)t_i}}{2}||\nabla{u_i}||_{2}^2-N e^{N t_i}\int_{\RR^N}G(u_i)dx \to 0.
\end{eqnarray}
\end{lemma}
\begin{proof}
Since $I_e(t,v)=I_e(t,-v)$ for all $(t,v)\in \RR \times H$, the group $K:=\{ (id,id), (id,-id)\}$ acts on the Banach space $\RR \times H$ isometrically. Due to the symmetry of $D_n$ and definition of $\Lambda_n^H$,
for any $\lambda(x)=(\tau(x), \gamma(x))\in \Lambda_n^H$, $(\tau(x), -\gamma(x))=(\tau(-x), \gamma(-x)) =\lambda(-x)$, which implies that the set $\lambda(D_n)$ is invariant under the group action $K$.
By the definition of $c_n^H$, there exists a sequence of $\{\gamma_i\}$ such that
$$ c_n^H \le \max_{x\in D_n} I(\gamma_i(x)) \le c_n^H +\frac{1}{i}.$$
Let $\lambda_i(x)= (0, \gamma_i(x))$, we see that $\lambda_i\in \Lambda_n^H$. Because $c_n^H=b_n^H$ and \eqref{relation}, we have
$$ b_n^H \le \max_{x\in D_n} I_e(\lambda_i(x)) \le b_n^H +\frac{1}{i}.$$
With the equivariant deformation lemma 3.1 and lemma 2.3 in Willem \cite{MW}, we obtain a sequence $\{(t_i, u_i)\}$ such that $I_e(t_i,u_i) \in [b_n^H-\frac{2}{i}, b_n^H+\frac{2}{i}]$ and
\begin{equation*}
| ( I_e)_t'(t,v)|+ ||(I_e)_v'(t,v))||_{H^*}=||I_e'(t,v)||_{\RR \times H^*}\le \frac{8}{\sqrt{i}}, \ \ \textrm{dist}_{\RR \times E_r}((t_i, u_i), \lambda_i(D_n))\le \frac{2}{\sqrt{i}}.
\end{equation*}
Immediately we have, as $i \to \infty$,  $I_e(t_i,u_i) \to c_n^H$ and $(I_e)_v'(t_i,u_i) \to 0$ strongly in $ H^{*}$. In addition we also have
$$ \frac{(N-2)e^{(N-2)t_i}}{2}||\nabla{u_i}||_{2}^2-N e^{N t_i}\int_{\RR^N}G(u_i)dx \to 0.  $$
Since $\lambda_i(x)= (0, \gamma(x))$, we find that
$\textrm{dist}_{\RR \times E_r}((t_i, u_i), \lambda_i(D_n))\ge |t_i|$. So $|t_i|\le \frac{2}{\sqrt{i}}\to 0$ as $i \to \infty$.
\end{proof}
Since $\lim_{i\to \infty}t_i= 0$, $\{u_i\}\subset H $ will become a Palais-Smale sequence for the minimax value $c_n^H$ if we show uniform boundedness of  $\{u_i\}$ in $H$.
Recall $g_{1}(s)=\max\{g(s)+ms^q, 0\}$ if $s\ge 0$ and
$g_{1}(s)=-g_{1}(-s)$ if $s\le 0$, and $G_1(s)=\int_{0}^{s}g_1(t)dt$. We also set $g_2(s):= g_1(s)-g(s)$ and $G_2(s):=\int_{0}^{s}g_2(t)dt$.
\begin{lemma}\label{bounded-1}
Let $\{(t_i,u_i)\}$ be the sequence from \eqref{Palais-Smale-1}. Then $\{u_i\}$ is uniformly bounded in $E$.
\end{lemma}
\begin{proof}
From \eqref{Palais-Smale-1}, by direct computation we have
\begin{equation}\label{bounded-g}
||\nabla{u_i}||_2^2 \to N c_n^H\ \ \  \ \textrm{and} \ \ \ \int_{\RR^N}G(u_i)dx \to \frac{N-2}{2} c_n^H.
\end{equation}
We divide the proof into the case $N\ge 3 $ and the case $N=2$. First we deal with the case $N\ge 3 $. Notice that there exists small $\delta>0$ such that $G_1(s)=0$ if $|s| \le \delta$, it follows that $ G_1(s) \le C_1 |s|^{2^*} $ for some $C_1>0 $. From \eqref{bounded-g}, for some $C_2>0$, we get
\begin{equation*}
\frac{m||u_i||_{q+1}^{q+1}}{q+1}\le \frac{N-2}{2} c_n^H+ \int_{\RR^N}G_2(u_i)dx+o(1)=\int_{\RR^N}G_1(u_i)dx\le C_1 ||u_i||_{2^*}^{2^{*}}\le C_2 ||\nabla{u_i}||_2^{2^*},
\end{equation*}
which implies that $||u_i||_{q+1}$ is uniformly bounded. Therefore $\{u_i\}$ is uniformly bounded in $E$ for $N\ge 3$.

Next we deal with the case $N=2$. According to our setting, we have $H=E_r$. From \eqref{Palais-Smale-1}, we have $\int_{\RR^2}G(u_i)dx \to 0$ and  $||\nabla{u_i}||_2^2 \to 2c_n^H$.
We need to show that $||u_i||_{q+1}$ is uniformly bounded, and this will be done by contradiction and rescaling arguments. Now we assume otherwise, $||u_i||_{q+1} \to \infty$ as $i \to \infty$.
Let $\delta_i=||(I_e)_v'(t,u_i))||_{E_r^*}$ and $w_i(y)=u_i(y/l_i)$, where $l_i=||u_i(x)||_{q+1}^{-(q+1)/2}$. In this way, we have, as $i\to \infty$,
\begin{equation}\label{bounded}
\delta_i \to 0,\ \ l_i \to 0, \ \ ||w_i(y)||_{q+1}=1, \ \ ||\nabla{w_i}||_2^2 \to 2c_n,\ \ \int_{\RR^2}G(w_i(y))dy \to 0.
\end{equation}
Since $\{ w_i\} $ is uniformly bounded in $E_r$, we may assume $w_i \rightharpoonup w_0$ weakly in $ E_r$ and pointwise a.e. in $\RR^2$.
\begin{claim}
$w_0$ is not a zero function in $E_r$.
\end{claim}
Indeed by definition, $g_2(s)\ge ms^q $ for $s\ge 0$, which means that $G_{2}(s)\ge \frac{m}{q+1}|s|^{q+1}$ for any $s\in\RR$. So we have
\begin{equation}\label{notzero-1}
\int_{\RR^2}G_{2}(w_i)dy\ge \frac{m}{q+1}\int_{\RR^2}|w_i|^{q+1}dy= \frac{m}{q+1}.
\end{equation}
Since $G(s)=G_1(s)-G_{2}(s)$ for any $s\in \RR$, by theorem~\ref{compactness-1}, we have
\begin{equation}\label{notzero-2}
\int_{\RR^2}G(w_i)dy+\frac{m}{q+1}\le \int_{\RR^2}G(w_i)dy+\int_{\RR^2}G_{2}(w_i)dy=\int_{\RR^2}G_{1}(w_i)dy \to \int_{\RR^2}G_{1}(w_0)dy.
\end{equation}
But $\int_{\RR^2}G(w_i)dy \to 0$, so we get $\int_{\RR^2}G_{1}(w_0)dy\ge \frac{m}{q+1}$, which implies $w_0$ is not zero.
\begin{claim}
$w_0$ is a zero function in $E_r$.
\end{claim}
Notice that $(I_e)_v'(t_i,u_i) \to 0$ strongly in $E_r^{*}$, for any $h(x)\in E_r$, we have
\begin{equation}\label{zero-1}
| \int_{\RR^2}\nabla{u_i} \nabla{h} dx - e^{2t_i}\int_{\RR^2}g(u_i) h dx |\le \delta_i(||\nabla{h(x)}||_2+||h(x)||_{q+1}).
\end{equation}
Now take $v(y)\in E_r$ and set $h(x)=v(y)$ with $y=l_ix$, we get
 $$\nabla_y{v}=\frac{1}{l_i}\nabla_x{h}, \ \ ||\nabla{h(x)}||_2=||\nabla{v(y)}||_2\ \textrm{ and } \ \ ||h(x)||_{q+1}=l_i^{-2/(q+1)}|| v(y)||_{q+1} . $$
Put this $h(x)$ into \eqref{zero-1} and make the change of variable from $x$ to $y$, and we have
\begin{equation}\label{zero-2}
| \int_{\RR^2}\nabla{w_i} \nabla{v} dy - \frac{e^{2t_i}}{l_i^2}\int_{\RR^2}g(w_i) v(y) dy |\le \delta_i(||\nabla{v(y)}||_2+l_i^{-2/(q+1)}|| v(y)||_{q+1}).
\end{equation}
By \eqref{inequality-4}, we have that $|g(w_i)|^{\frac{q+1}{q}}\le C(|w_i|^{q+1}+|w_i|^{q+1}(e^{\frac{q+1}{q}\alpha w_i^2}-1))$. By \eqref{bounded}, $w_i$ is uniformly bounded in $E$.
Through theorem~\ref{TRU-MOS-F}, choosing $\alpha$ small, we see that $g(w_i)$ is uniformly bounded in $L^{(q+1)/q}(\RR^2)$.
Observing that $g(w_i)$ converges to $g(w_0)$ pointwise a.e. in $\RR^2$ and $v(y)\in L^{q+1}(\RR^2)$, by weak convergence we have, as $i \to \infty$,
\begin{equation}\label{zero-3}
\int_{\RR^2}g(w_i) v(y) dy \to \int_{\RR^2}g(w_0) v(y) dy.
\end{equation}
Multiplying both sides of \eqref{zero-2} by $l_i^2$, we find
\begin{equation*}
| l_i^2\int_{\RR^2}\nabla{w_i} \nabla{v} dy - e^{2t_i}\int_{\RR^2}g(w_i) v(y) dy |\le \delta_i(l_i^2||\nabla{v(y)}||_2+l_i^{2q/(q+1)}|| v(y)||_{q+1}).
\end{equation*}
Taking the limit as $i \to \infty$, we have $\int_{\RR^2}g(w_0) v(y) dy=0$, which means $g(w_0)=0$. By (g1), we find $w_0=0$, which is a contradiction with the previous claim. Therefore $\{u_i\}$ is uniformly bounded in $E_r$
\end{proof}
\begin{remark}\label{Palais-Smale-2}
As we pointed out before lemma~\ref{bounded-1} that $\{u_i\}$ became a Palais-Samle sequence after we proved the uniform boundedness of  $\{u_i\}$. Indeed from \eqref{Palais-Smale-1}, we have
$I_e(t_i,u_i) \to c_n^H$, $(I_e)_v'(t_i,u_i) \to 0$ strongly in $ H^{*}$. Since $u_i$ is uniformly bounded, for any $v \in H$, we have
\begin{equation*}
|( e^{(N-2)t_i}-1)\int_{\RR^N}\nabla{u_i} \nabla{v} dx |=o(1)||\nabla{v}||_2, \ \ |(e^{Nt_i}-1)\int_{\RR^N}g(u_i)v dx |= o(1)(||v||_{q+1} +||\nabla{v}||_2),
\end{equation*}
which implies $I(u_i) \to c_n^H$ and $||I'(u_i)||_{H^*} \to 0$ as $i \to \infty$.
Moreover this Palais-Smale sequence has an additional relation
\begin{equation}\label{Pohozaev-1}
P(u_i)=\frac{N-2}{2}||\nabla{u_i}||_{2}^2-N \int_{\RR^N}G(u_i)dx \to 0.
\end{equation}
\end{remark}
Now we deal with the convergence of this Palais-Samle sequence $\{u_i\}\subset H$. We first treat the case when $H=X_\tau \cap E_1$. Let $x=(x^1,x^2,x^3)\in \RR^N=\RR^M \times \RR^M \times \RR^{N-2M}$.
\begin{theorem}\label{convergence-weak}
Let $H=X_\tau \cap E_1$ and  $\{(t_i,u_i)\}$ be the sequence from \eqref{Palais-Smale-1}. If $N\ge 5$ and $N-2M\ge 1$, there exists $ \{y_i=(0,0, y_i^3)\} \subset \RR^N $, such that $u_i(\cdot+y_i)$ converges weakly to some nontrivial function ${u}_n^H \in X_{\tau}\cap E_1$, which is a solution to \eqref{main} with $I({u}_n^H)\le c_n^H$.
\end{theorem}
\begin{proof}
We claim that there exists $r_0>0$ such that
\begin{equation}\label{notzero-1}
\lim_{i \to \infty} \sup_{x^3\in \RR^{N-2M}}\int_{B((0,0,x^3),r_0)}|u_i|^{q+1} dx >0.
\end{equation}
Otherwise for all $r>0$, we have $\lim_{i \to \infty} \sup_{x^3\in \RR^{N-2M}}\int_{B((0,0,x^3),r)}|u_i|^{q+1} dx =0$, which in turn implies that
\begin{equation}\label{TM-1}
\lim_{i \to \infty} \sup_{x \in \RR^{N}}\int_{B(x ,1)}|u_i|^{q+1} dx =0.
\end{equation}
Indeed if it is not true, then there exists a sequence $\{ x_i=(x_i^1, x_i^2, x_i^3)\} \subset \RR^N$ such that $\int_{B(x_i ,1)}|u_i|^{q+1} dx \ge C $ for some $C>0$. Because $v_i \in E_1$ and $||u_i||_{q+1}$ is uniformly bounded, there is $r_1>0$ such that $|x_i^1|^2+|x_i^2|^2 \le r_1^2$. It follows that
\begin{equation*}
\int_{B((0,0,x^3_i),r_1+1)}|u_i|^{q+1} dx \ge \int_{B(x_i ,1)}|u_i|^{q+1} dx \ge C,
\end{equation*}
which is impossible. So \eqref{TM-1} is true. By lemma I.1. from \cite{LION}, we have $\lim_{i \to \infty} ||u_i||_p =0$ for $q+1<p<2^*$. Repeating the proof of the second part of theorem~\ref{compactness-1}, we get $ \int_{\RR^N} G_1(u_i) dx  \to 0$ as $i \to \infty$. From \eqref{Pohozaev-1}, we find that
\begin{equation*}
\frac{N-2}{2}||\nabla{u_i}||_{2}^2\le \frac{N-2}{2}||\nabla{u_i}||_{2}^2+N \int_{\RR^N}G_2(u_i)dx=N \int_{\RR^N}G_1(u_i)dx+o(1)=o(1).
\end{equation*}
In viewing of \eqref{bounded-g}, we reach a contradiction. So there exists $r_0>0$ such that \eqref{notzero-1} is true. It follows that there exists $ \{y_i=(0,0, y_i^3)\} \subset \RR^N $, such that $\tilde{u}_i: =u_i(\cdot+y_i)$ converges weakly and point-wise almost everywhere in $\RR^N$ to some nontrivial function ${u}_n^H\in X_{\tau}\cap E_1$. Moreover $\{\tilde{u}_i \}\subset H $ is also a Palais-Samle sequence, i.e.  as $i \to \infty$, $||I'(\tilde{u}_i)||_{H^*} \to 0$ and
\begin{equation*}
\frac{1}{2}||\nabla{\tilde{u}_i}||_{2}^2- \int_{\RR^N}G(\tilde{u}_i)dx =I(\tilde{u}_i) \to c_n^H, \ P(\tilde{u}_i)=\frac{N-2}{2}||\nabla{\tilde{u}_i}||_{2}^2-N \int_{\RR^N}G(\tilde{u}_i)dx \to 0.
\end{equation*}
Then we have $N c_n^H+o(1)=||\nabla{\tilde{u}_i}||_{2}^2$, which shows that $N c_n^H \ge ||\nabla{{u}_n^H}||_{2}^2$. Finally by weak convergence we have $I'({u}_n^H)v =0$ for any $v\in H$. Indeed we only need to show that
\begin{equation}\label{TM-2}
\lim_{i \to \infty} \int_{\RR^N} |g(\tilde{u}_i) v| dx = \int_{\RR^N} |g({u}_n^H) v| dx .
\end{equation}By (g1) and (g2), we have $$|g(\tilde{u}_i)v|\le (C_2 |\tilde{u}_i|^{q} |v| +\epsilon |\tilde{u}_i|^{2^*-1}|v|).$$
By weak convergence and point-wise convergence, we have $\int_{\RR^N} |\tilde{u}_i|^{q} |v|dx \to \int_{\RR^N} |{u}_n^H|^{q} |v|dx$. Repeating the proof of the second part of theorem~\ref{compactness-1}, we have
$\limsup_{i \to \infty} \int_{\RR^N} |g(\tilde{u}_i) v| dx \le \int_{\RR^N} |g({u}_n^H) v| dx $. Combining this with Fatou's lemma, we get \eqref{TM-2}, which in turn implies
\begin{equation}\label{TM-3}
\int_{\RR^N} g(\tilde{u}_i) v dx \to \int_{\RR^N} g({u}_n^H) v dx .
\end{equation}
 Therefore we have $I'({u}_n^H)v =0$ for any $v\in H$. By principle of symmetric criticality \cite{PAL}, ${u}_n^H$ is a solution to \eqref{main}. From Pohozaev identity, we have $I({u}_n^H)= \frac{1}{N}||\nabla{{u}_n^H}||_{2}^2\le c_n^H$.
\end{proof}

\begin{theorem}\label{compactness-2}
Assume $H=E_2\cap X_\tau$ with $N=4$ or $N\ge 6$, and $N-2M\neq 1$, or $H=E_r$. Let $\{(t_i,u_i)\}$ be the sequence from \eqref{Palais-Smale-1}. Then $u_i$ converges strongly to some nontrivial function ${u}_n^H \in H$, which is a solution to \eqref{main} with $I({u}_n^H)= c_n^H$.
\end{theorem}
\begin{proof}
By lemma~\ref{bounded-1}, taking a subsequence if necessary, we assume $u_i$ converges to $u_n^H$ weakly in $H$ and pointwise a.e. in $\RR^N$ as $i \to \infty$.
For any $v \in H$, similar to \eqref{zero-3}, \eqref{TM-2} and \eqref{TM-3}, we have
\begin{equation*}
\lim_{i\to \infty} \int_{\RR^N}g(u_i) vdx = \int_{\RR^N}g(u_n^H) v dx,
\end{equation*}
which implies, by weak convergence, $\int_{\RR^N} \nabla{u_n^H} \nabla{v} dx =\int_{\RR^N}g(u_n^H)vdx $. In particular
\begin{equation}\label{zero-4}
||\nabla{u_n^H}||_2^2=\int_{\RR^N} \nabla{u_n^H} \nabla{u_n^H} dx =\int_{\RR^N}g(u_n^H)u_n^Hdx.
\end{equation}
By theorem~\ref{compactness-1}, we have
\begin{equation}\label{zero-5}
\lim_{i\to \infty} \int_{\RR^N}g_1(u_i) u_i dx = \int_{\RR^N}g_1(u_n^H) u_n^H dx.
\end{equation}
Since $|I'(u_i)u_i| \to 0$ as $i\to \infty$, we have, by \eqref{zero-4} and \eqref{zero-5},
\begin{eqnarray}
& &\limsup_{i\to \infty}(||\nabla{u_i}||_2^{2}+\int_{\RR^N}g_2(u_i)u_idx) =\limsup_{i\to \infty}(||\nabla{u_i}||_2^{2}+\int_{\RR^N}g_2(u_i)u_idx-I'(u_i)u_i) \nonumber \\
& & =\limsup_{i\to \infty}\int_{\RR^N}g_1(u_i)u_idx=\int_{\RR^N}g_1(u_n^H) u_n^H dx=||\nabla{u_n^H}||_2^2+\int_{\RR^N}g_2(u_n^H)u_n^Hdx   \nonumber
\end{eqnarray}
Notice that $u_i \to u_n^H$ pointwise a.e. in $\RR^N$ and $g_2(s)s \ge 0 $ for all $s \in \RR$, by weak lower semi-continuity and Fatou's lemma, we get
\begin{eqnarray}
& & ||\nabla{u_n^H}||_2^{2}+\int_{\RR^N}g_2(u_n^H)u_n^Hdx \le \liminf_{i\to \infty}||\nabla{u_i}||_2^{2}+\liminf_{i\to \infty}\int_{\RR^N}g_2(u_i)u_idx \nonumber \\
& & \le \liminf_{i\to \infty}(||\nabla{u_i}||_2^{2}+\int_{\RR^N}g_2(u_i)u_idx) \le \limsup_{i\to \infty}(||\nabla{u_i}||_2^{2}+\int_{\RR^N}g_2(u_i)u_idx).  \nonumber
\end{eqnarray}
These two calculations above show that we must have $\lim_{i\to \infty}||\nabla{u_i}||_2=||\nabla{u_n^H}||_2$ and $\lim_{i\to \infty}\int_{\RR^N}g_2(u_i)u_idx=\int_{\RR^N}g_2(u_n^H)u_n^Hdx$.
Because $g_2(s)s\ge m|s|^{q+1} $ for $s\in \RR$ and $u_i \to u_n^H$ pointwise a.e. in $\RR^N$, by Lebesgue convergence theorem, we have $\lim_{i\to \infty}||u_i||_{q+1}=||u_n^H||_{q+1}$.
Therefore we have $\lim_{i\to \infty}||u_i||_{E}=||u_n^H||_{E}$. Since $u_i$ converges to $u_n^H$ weakly in $H$, $u_i$ converges to $u_n^H$ strongly in $H$. Moreover $I(u_n^H)=c_n^H$.
\end{proof}

\section{Analysis of the energy level}
In this section we are going to analyze the energy level $c_1^H$ when $H=E_r$, or $H=X_\tau \cap E_1$ with $N\ge 5$ and $N-2M\ge 1$, or $H=E_2\cap X_\tau$ with $N=4$ or $N\ge 6$, and $N-2M\neq 1$. From \eqref{minimax-1},
noticing that $I(v)=I(-v)$, $D_1=[-1, 1]$ and $\gamma(0)=0$ for any $\gamma \in \Gamma_1^H$, we have
\begin{equation}\label{curve-3}
c_{1}^H=\inf_{\gamma \in \Gamma_{1}^H} \sup_{x\in [0, 1]} I(\gamma(x))=\inf_{\gamma \in \Gamma_H} \sup_{x\in [0, 1]} I(\gamma(x)),
\end{equation}
where $\Gamma_H=\{ \gamma\in C([0,1], H)\ |\ \gamma(0)=0 \ ; \ I(\gamma(1))<0 \}$. It is easy to see that $|\gamma| \in \Gamma_H$ if $\gamma \in \Gamma_H$.
\begin{lemma}\label{curve-1}
Let $H=E_r$. Then there exists $\bar{u}_1 \in E_r$ such that $\bar{u}_1 \ge 0$, $I'(\bar{u}_1 )=0$ and $I(\bar{u}_1 )=c_1^{E_r}$. In particular $\bar{u}_1 $ is a nonnegative solution to \eqref{main}.
\end{lemma}
\begin{proof}
Basically we will repeat what we did in section 3. We define $b_H$:
\begin{equation*}
b_H=\inf_{\lambda \in \Lambda_H} \sup_{x\in [0,1]} I_e(\lambda(x)),
\end{equation*}
where $\Lambda_H=\{ \lambda=(\tau(x), \gamma(x)) \in C([0,1], \RR \times E_r)\ | \ \gamma(0)=0;  \ I_e(\tau(1), \gamma(1))<0 \}$. Repeat lemma~\ref{coincide-1}, we have $\Lambda_H \neq \emptyset $ and $c_1^H =b_H$.
Choose a sequence of $\{\gamma_i\}\subset \Gamma_H$ such that
$$  c_1^H \le \max_{x\in [0,1]} I(\gamma_i(x)) \le c_1^H +\frac{1}{i}, \ \ \gamma_i(x)\ge 0 \ \textrm{ for any $x\in [0,1]$}.$$
Let $\lambda_i(x)\triangleq (0, \gamma_i(x))$, we see that $\lambda_i\in \Lambda_H$. Because $c_1^H=b_H$ and \eqref{relation}, we have
$$ c_1^H \le \max_{x\in [0,1]} I_e(\lambda_i(x)) \le c_1^H +\frac{1}{i}.$$
Following lemma~\ref{coincide-2}, theorem~\ref{bounded-1} and lemma~\ref{compactness-2}, we have a sequence $\{(t_i, v_i) \}\subset \RR \times E_r$ such that, as $i \to \infty$, $t_i\to 0$,
$\textrm{dist}_{\RR \times E_r}((t_i, v_i), \lambda_i([0,1])) \to 0$,
$I'(v_i) \to 0 $ strongly in $E_r^*$ and $I(v_i) \to c_1^H$. Moreover $\{v_i\}$ is bounded and compact in $E_r$. Now we assume $ v_i \to \bar{u}_1$ strongly in $E_r$, so we have $I'(\bar{u}_1)=0$ and $I(\bar{u}_1)=c_1^H$.
Notice that $\gamma_{i}(x) \ge 0$ for any $x\in [0,1]$ and $\textrm{dist}_{\RR \times E_r}((t_i, v_i), \lambda_i([0,1])\ge \int_{\RR^N}( \max\{0,-v_i\})^{q+1}dx$, so we have, as $i \to \infty$,
 $\int_{\RR^N} (\max\{0,-v_i\})^{q+1}dx \to 0$. We see that  $\int_{\RR^N} (\max\{0,-\bar{u}_1\})^{q+1}dx =0$, which implies $\bar{u}_1$ is nonnegative.
\end{proof}
\begin{remark}\label{positive-1}
If $q\ge 1 $, the strong maximum principle applies and any nontrivial nonnegative solution of \eqref{main} in $E$ is strictly positive in $\RR^N$. But if $0<q<1$, the strong maximum principle does not applies. In fact, as proved in \cite{SQ5}, by comparison principle, any solution of \eqref{main} in $E$ has compact support.
\end{remark}
Let $H=E_r$, we want to construct a continuous curve in $\Gamma_H$ to show that $c_1^{E_r}$ is the least energy level for all nontrivial solutions. For the case $q=1$, the construction of such curve can also be found in \cite{JT}. If $0<q<1$ and $N=2$, the nonlinear term $g(s)$ is not Lipschitz at $s=0$, the construction in \cite{JT} is not directly applicable and we have to do some adjustments.
\begin{lemma}\label{curve-2}
Let $w$ be any nontrivial solution  of \eqref{main} in $E$ and $w^*$ be the Schwarz spherical rearrangement of $|w|$. Then there exist two continuous curves $\gamma(t):[0,1] \to E$ and $\gamma^*(t): [0,1] \to E_r$ such that
$ \gamma \in \Gamma_H$, $ \gamma^* \in \Gamma_H$, $\gamma(\bar{t}) =w$ and $\gamma^*(\bar{t}) =w^*$ for some $\bar{t}\in (0,1)$. Moreover
\begin{equation*}
\ I(\gamma^*(t))\le I(\gamma(t)) \ \textrm{for any } \ t \in [0,1] \ \textrm{and}\ \max_{t \in [0,1]} I(\gamma^*(t))\le \max_{t \in [0,1]} I(\gamma(t))=I(w).
\end{equation*}
\end{lemma}
\begin{proof}
Sine $w^*$ is the Schwarz spherical rearrangement of $|w|$, we have
\begin{equation}\label{rearrangement}
||\nabla{w^*}||_2\le ||\nabla{w}||_2 \ \  \textrm{and} \ \ \int_{\RR^N}G(w^*)dx =\int_{\RR^N}G(|w|)dx =\int_{\RR^N}G(w)dx.
\end{equation}
And we also have $\int_{\RR^N}G(\tau w^*)dx =\int_{\RR^N}G(\tau w)dx$ for any $\tau>0$. We first construct the curves for the case $N\ge 3$. Let us define $\gamma(t):=w(x/t) $ for $t>0$ and $\gamma(t):=0$, and define $\gamma^*(t):=w^*(x/t) $ for $t>0$ and $\gamma^*(t):=0$. Because $w$ is a solution of \eqref{main}, we find
\begin{equation*}
P(w)=\frac{N-2}{2}||\nabla{w}||_{2}^2-N \int_{\RR^N}G(w)dx =0.
\end{equation*}
It follows that
\begin{equation*}
I(\gamma^*(t))\le I(\gamma(t))=\frac{t^{N-2}}{2}||\nabla{w}||_2^2-t^N\int_{\RR^N}G(w)dx=\frac{||\nabla{w}||_2^2}{2}(t^{N-2}-\frac{N-2}{N}t^N),
\end{equation*}
which means that $I(\gamma(L))< 0$ for some large $L>1$ and  $\max_{t \in [0,L]} I(\gamma(t))=I(w)$. By rescaling we are done for the case $N\ge 3$.

Now we deal with the case $N=2$. In this case the Pohozaev identity \eqref{Pohozaev-2} reduces to $\int_{\RR^N}G(w)dx=0$.  We divide the construction into four steps.
Set $C_1:=\max_{\tau\in [0,1]} \int_{\RR^2} |G(\tau w)| +|g(\tau w)\tau w| dx  $ and $t_1:=\min \{ ||\nabla{w}||_2/(2\sqrt{C_1}), \ 1/2\}$.\\
\textbf{Step one:} Let $\gamma_1(\tau):= \tau w(x/(\tau t_1)): (0,1] \to E $ and $\gamma_1^*(\tau)=: \tau w^*(x/(\tau t_1)): (0,1] \to E_r $. Define $\gamma(0)=\gamma^*(0)=0$.
Clearly $\gamma_1(\tau)$ and $\gamma_1^*(\tau)$ are continuous, and we also have
\begin{equation*}
I(\gamma_1(\tau))=\frac{\tau^2}{2}||\nabla{w}||_2^2-\tau^2t_1^2\int_{\RR^2}G(\tau w)dx=\frac{\tau^2}{2}||\nabla{w}||_2^2-\tau^2t_1^2\int_{\RR^2}G(\tau w^*)dx\ge I(\gamma_1^*(\tau)).
\end{equation*}
By choice of $t_1$ we find $\frac{d(I(\gamma_1(\tau)))}{d\tau}=\tau( ||\nabla{w}||^2-t_1^2\int_{\RR^2} 2 G(\tau w)+g(\tau w)\tau w dx ) \ge 0$ for all $\tau \in [0,1]$.
Therefore we see that $I(\gamma_1(\tau))$ is increasing and $(1/2)||\nabla{w}||_2^2=I(\gamma_1(1))\ge I(\gamma_1(\tau))\ge I(\gamma_1^*(\tau))$ for all $\tau \in [0,1]$.
\\
\textbf{Step two:} Let $\gamma_2(t) := w(x/t): [t_1, 2] \to E $ and $\gamma_2^*(t):= w^*(x/t): [t_1, 2] \to E_r $. Notice that $t_1\le 1/2$,
$\gamma_2(1)=w$, and both $\gamma_2(t)$ and $\gamma_2^*(t)$ are well-defined and continuous.
By direct calculation, we have $I(\gamma_2(t))=(1/2)||\nabla{w}||_2^2\ge (1/2)||\nabla{w^*}||_2^2=I(\gamma_2^*(t))$. Set $M(s)=I(s\gamma_2(2))$, we have $M(s)=\frac{s^2}{2}||\nabla{w}||_2^2-4\int_{\RR^2}G(s w)dx$ and
\begin{equation*}
M'(1)=||\nabla{w}||_2^2-4\int_{\RR^2}g(w)wdx=-3||\nabla{w}||_2^2<0,
\end{equation*}
which implies $M(s)$ is strictly decreasing in $[1, s_1]$ for some $s_1>1$ close to $1$. Since $M(1)>M(s_1)$, we also have $\int_{\RR^2}G(s_1w)dx >0$.
\\
\textbf{Step three:} Let $\gamma_3(s):= s\gamma_2(2)=sw(x/2): [1, s_1] \to E $ and $\gamma_3^*(t):= sw^*(x/2): [1, s_1] \to E_r $.
It is clear that $\gamma_3(s)$ and $\gamma_3^*(s)$ are continuous, and $I(\gamma_3(s))$ is decreasing. Moreover for any $s \in [1, s_1]$, we have
\begin{equation*}
(1/2)||\nabla{w}||_2^2=I(\gamma_3(1))\ge I(\gamma_3(s))= \frac{s^2}{2}||\nabla{w}||_2^2-4\int_{\RR^2}G(s w)dx \ge I(\gamma_3^*(s)).
\end{equation*}
\\
\textbf{Step four:} Let $\gamma_4(t):= s_1 w(x/t): [2, t_2] \to E $ and $\gamma_4^*(t):= s_1w^*(x/t): [2, t_2] \to E_r $, where $t_2:=\frac{s_1||\nabla{w}||_2}{\sqrt{\int_{\RR^2}G(s_1w)dx}}+2$.
By step two, we know that  $\int_{\RR^2}G(s_1w)dx >0$, which implies $\gamma_4(t)$ and $\gamma_4^*(t)$ are well-defined. By direct computation we see that $I(\gamma_4(t))$ is strictly decreasing and
\begin{equation*}
(1/2)||\nabla{w}||_2^2=I(\gamma_3(1))>I(\gamma_3(s_1))=I(\gamma_4(2))\ge I(\gamma_4(t))\ge I(\gamma_4^*(t)), \ I(\gamma_4(t_2))<0.
\end{equation*}
Finally connecting $\gamma_1$,$\gamma_2$,$\gamma_3$ and $\gamma_4$, by rescaling we can get the desired curve $\gamma$. Similarly We can connect $\gamma_1^*$,$\gamma_2^*$,$\gamma_3^*$ and $\gamma_4^*$ to get$\gamma^*$.
\end{proof}
\begin{remark}\label{weak-assumption}
From the construction above, we see that the assumption that $w$ is a solution to \eqref{main} can be weakened. Actually when $N\ge 3$, we only need $w$, satisfying Pohozaev identity $P(w)=0$, to finish the construction. When $N=2$, we need more. The function $w$ has to satisfy $P(w)=0$ and stay inside the Nehari manifold, i.e. $||\nabla{w}||_2^2=\int_{\RR^N} g(w)wdx $.
\end{remark}
\begin{theorem}\label{least-energy-solution}
Let $H=E_r$, or $H=X_\tau \cap E_1$ with $N\ge 5$ and $N-2M\ge 1$, or $H=E_2\cap X_\tau$ with $N=4$ or $N\ge 6$, and $N-2M\neq 1$. Then $c_1^H= \inf\{ \ I(v)\ | \ I'(v)=0;\ v\in H-\{0\} \} $ and
\begin{equation*}
c_1^{E_r}= \inf\{ \ I(v)\ | \ I'(v)=0;\ v\in E-\{0\} \} = \inf\{ \ I(v)\ | \ P(v)=0;\ v\in E-\{0\} \}.
\end{equation*}
\end{theorem}
\begin{proof}
By theorem~\ref{convergence-weak} and theorem~\ref{compactness-2}, there exists solution $u_1^H \in H$ such that $I(u_1^H)\le c_1^H$.
Then by the definition of $c_1^H$ from \eqref{curve-3} and lemma~\ref{curve-2}, we have  $$c_1^H \le \inf\{ \ I(v)\ | \ I'(v)=0;\ v\in H-\{0\} \} \le I(u_1^H)\le c_1^H. $$
It follows that $c_1^H= \inf\{ \ I(v)\ | \ I'(v)=0;\ v\in H-\{0\} \} $. Moreover by lemma~\ref{curve-2} again, we have
$$c_1^{E_r} \le \inf\{ \ I(v)\ | \ I'(v)=0;\ v\in H-\{0\} \} \le \inf\{ \ I(v)\ | \ I'(v)=0;\ v\in E_r \}=c_1^{E_r}, $$
that is $c_1^{E_r}= \inf\{ \ I(v)\ | \ I'(v)=0;\ v\in E-\{0\} \}$. So we are left to prove $c_1^{E_r}= \inf\{ \ I(v)\ | \ P(v)=0;\ v\in E-\{0\} \} $. Since any solution of \eqref{main} satisfies Pohozaev identity, we have
$c_1^{E_r}\ge \inf\{ \ I(v)\ | \ P(v)=0;\ v\in E-\{0\} \} $. For the case $ N\ge 3$, by remark~\ref{weak-assumption} and lemma~\ref{curve-2}, $c_1^{E_r}\le \inf\{ \ I(v)\ | \ P(v)=0;\ v\in E-\{0\} \} $. So we are done in this situation.

For the case $N=2$, following the arguments in \cite{BGK}, we show that the constrained minimization problem $\inf\{ \frac{1}{2} ||\nabla{v}||_2^2 \ | \ v\in S \ \}$, $S:= \{ \ v\in E-\{0\} \ | \ \int_{\RR^2} G(v) dx =0,  \}$ admit a solution to \eqref{main}. Indeed by lemma~\ref{odd-mapping}, there is a function $v_1\in E_r$ such that $\int_{\RR^2}G(v_1)dx>0$. From lemma~\ref{makesense}, by continuity there is $0<t<1$ such that $\int_{\RR^2}G(tv_1)dx=0$, which shows that $S$ is not empty. Set $K:=\inf\{ \frac{1}{2} ||\nabla{v}||_2^2 \ | \ v\in S \ \}$. From lemma~\ref{makesense}, we know $K>0$. Due to the Schwarz spherical rearrangement and rescaling properties, we may assume that the constrained minimization problem has a nonnegative minimizing sequence $\{ v_i \} \subset E_r$ such that as $i \to \infty$,
\begin{equation*}
||\nabla{v_i}||_2^2 \to 2K,\ ||v_i||_{q+1}=1,\ \int_{\RR^2} G(v_i)dx =0.
\end{equation*}
It follows that $v_i $ converges to $w $ weakly for some $w \in E_r$. By theorem~\ref{compactness-1}, we get as $i \to \infty$,
\begin{equation*}
\frac{m}{q+1} \le \frac{m}{q+1}\int_{\RR^2}v_i^{q+1}dx\le  \int_{\RR^2} G_2(v_i)dx =\int_{\RR^2} G_1(v_i)dx \to  \int_{\RR^2} G_1(w)dx,
\end{equation*}
which tells that $w$ is nonzero function in $E_r$. From weak convergence, Fatou's lemma and theorem~\ref{compactness-1}, we have $0<||\nabla{w}||_2^2 \le 2K$ and $\int_{\RR^2} G(w)dx \ge 0$.
We can rule out $\int_{\RR^2} G(w)dx >0 $. Indeed from lemma~\ref{makesense}, by continuity there is $0<t<1$ such that $\int_{\RR^2}G(tw)dx=0$. Obviously $ ||\nabla{tw}||_2^2 < 2K $ and this is a contradiction.
Hence we have $\int_{\RR^2} G(w)dx =0 $, in turn we find $||\nabla{w}||_2^2 = 2K$. By Lagrange multiplier rule, there exists $\lambda \in \RR$ such that $ \int_{\RR^2}\nabla{w} \nabla{v} dx =\lambda \int_{\RR^2}g(w)v dx $. Following the arguments in \cite{BGK} or in \cite{SQ5}, we can show that $\lambda <0$ is impossible and $\lambda$ has to be positive. Finally we see that $w(x/\sqrt{\lambda})$ is a solution to \eqref{main}. It is clear $P(w(x/\sqrt{\lambda}))=\int_{\RR^2} G(w(x/\sqrt{\lambda}))dx=0$ and $I(w(x/\sqrt{\lambda}))=\inf\{ \ I(v)\ | \ P(v)=0;\ v\in E-\{0\} \} $. So we have
\begin{equation*}
c_1^{E_r} \le \inf\{ \ I(v)\ | \ P(v)=0;\ v\in E-\{0\} \} \le \inf\{ \ I(v)\ | \ I'(v)=0;\ v\in E-\{0\} \}=c_1^{E_r},
\end{equation*}
that is $c_1^{E_r}= \inf\{ \ I(v)\ | \ I'(v)=0;\ v\in E-\{0\} \} = \inf\{ \ I(v)\ | \ P(v)=0;\ v\in E-\{0\} \}$.
\end{proof}
\begin{remark}
Let $H=E_r$. From theorem~\ref{least-energy-solution}, we know that $c_1^{E_r}$ is the least energy level for solution to \eqref{main}. By lemma~\ref{curve-1}, \eqref{main} has one nonnegative least energy solution $\bar{u}_1 \in E_r$. In fact from \cite{BJM} we know that each least energy solution is radial up to translation and each least energy solution has constant sign. Since $g$ is an odd function, we can pick any least energy solution to be nonnegative.
\end{remark}
Now we give the proof of theorem~\ref{existence}.
\begin{proof}[Proof of theorem~\ref{existence}]
In this case we know $H=E_r$. Theorem~\ref{compactness-2} produces a sequence of solutions $\{ u_n^H \} $ to \eqref{main} with $I(u_n^H )=c_n^H$. Theorem~\ref{least-energy-solution} shows that $c_1^{E_r}$ is the least energy level. Set $\bar{u}_n:=u_n^H$. From lemma~\ref{curve-1}, we can choose $\bar{u}_1$ to be nonnegative. Theorem~\ref{unbounded-2} in the next section shows that $c_n^{E_r} \to \infty$ as $n \to \infty$.
\end{proof}
Next we give the proof of theorem~\ref{nonradial-existence}.
\begin{proof}[Proof of theorem~\ref{nonradial-existence}]
(i) In this case we know $H=X_{\tau}\cap E_1$. Take $n=1$ in Theorem~\ref{convergence-weak}, we find a solution ${u}_1^H \in H$ with $I({u}_1^H )\le c_1^H$.  Set $\tilde{u}:={u}_1^H$. By theorem~\ref{least-energy-solution}, we see that $\tilde{u}=c_1^H$, which minimizes the energy functional $I$ over all nontrivial solutions in $X_{\tau}\cap E_1$. We are left to prove $I(\tilde{u})>2 c_1^{E_r}$.
Knowing that $x=(x^1,x^2,x^3)\in \RR^N=\RR^M \times \RR^M \times \RR^{N-2M}$, we write
\begin{equation*}
\Omega^1=\{ x\in \RR^N \ | \ |x^1|>|x^2| \} \ \ \textrm{and}\ \ \Omega^2=\{ x\in \RR^N \ | \ |x^1|<|x^2| \}.
\end{equation*}
By standard regularity theory \cite{GT}, we know $\tilde{u}\in W_{loc}^{2,p}(\RR^N)$ for any $p>1$, so we have $\tilde{u}\in C_{loc}^{1,\alpha}(\RR^N)$ for some $0<\alpha<1$. Since $\tilde{u}\in X_{\tau}\cap E_1$, we have
$\tilde{u}(x)=0 $ if $|x^1|=|x^2|$. Let $\chi(\Omega^i)$ be the characteristic function of $\Omega^i$, $i=1,\ 2$. Then we have $\chi(\Omega^i)\tilde{u}\in E $, and both of them are nontrivial. By symmetry
\begin{equation*}
0=P(\tilde{u})=P(\chi(\Omega^1)\tilde{u})+P(\chi(\Omega^2)\tilde{u})=2P(\chi(\Omega^1)\tilde{u}),
\end{equation*}
which implies $P(\chi(\Omega^1)\tilde{u})=0$. In turn, by theorem~\ref{least-energy-solution}, we have
\begin{equation*}
I(\tilde{u})=I(\chi(\Omega^1)\tilde{u})+I(\chi(\Omega^2)\tilde{u})=2I(\chi(\Omega^1)\tilde{u}\ge 2 c_1^{E_r}.
\end{equation*}
It follows that $I(\tilde{u})\ge 2 c_1^{E_r}$. We rule out the possibility that $I(\tilde{u})= 2 c_1^{E_r}$. Actually if $I(\tilde{u})= 2 c_1^{E_r}$, then $I(\chi(\Omega^1)\tilde{u})=c_1^{E_r}$. By lemma~\ref{curve-2} and remark~\ref{weak-assumption}, for some $0<t_0<1$, we can construct a curve $\gamma(t): [0,1]\to E$ such that $\gamma(t_0)=\chi(\Omega^1)\tilde{u}$ and $c_1^{E_r}=I(\gamma(t_0))>I(\gamma(t))$ for $t\neq t_0$. We conclude that $\chi(\Omega^1)\tilde{u}$ must be one least energy solution to \eqref{main}, otherwise by lemma 2.3 in \cite{MW}, we can decrease the least energy level $c_1^{E_r}$. If $q\ge 1$, from \cite{BJM} and strong maximum principle, this least energy solution $\chi(\Omega^1)\tilde{u}$ must be strictly positive or negative over the whole space $\RR^N$, this is impossible. If $0<q<1$,  from \cite{BJM} and the symmetry of $\tilde{u}$, this least energy solution $\chi(\Omega^1)\tilde{u}$ is radially symmetric about some point $(0,0, y^3_1) \in \RR^N$ and least energy solution $\chi(\Omega^2)\tilde{u}$ is radially symmetric about some point $(0,0, y^3_2) \in \RR^N$. Both of them have constant sign and have $\mathcal{O}(2M)$ symmetry, which means that $\tilde{u}=\chi(\Omega^1)\tilde{u}+\chi(\Omega^2)\tilde{u} \notin X_{\tau}$. We also reach a contradiction. Therefore $I(\tilde{u})> 2 c_1^{E_r}$.

(ii) In this case we have $H=X_{\tau}\cap E_2$. According to theorem~\ref{compactness-2}, there exists a sequence of solutions $\{ u_n^H \} $ to \eqref{main} with $I(u_n^H )=c_n^H$.  Set $\hat{u}_n:=u_n^H$. By theorem~\ref{least-energy-solution}, we see that $\hat{u}_1$ minimizes the energy functional $I$ over all nontrivial solutions in $X_{\tau}\cap E_2$. Modifying the proof in (i) accordingly, we can show that $I(\hat{u}_1)> 2 c_1^{E_r}$.  Theorem~\ref{unbounded-2} in the next section shows that $c_n^{H} \to \infty$ as $n \to \infty$.
\end{proof}

\section{Comparison functional and unboundedness of minimax values}
Let $H=E_r$, or $H=E_2\cap X_\tau$ with $N=4$ or $N\ge 6$, and $N-2M\neq 1$. In this section we are going to show that  $\lim_{n \to\infty}c_n^H=\infty$, where $c_n^H$ is defined in \eqref{minimax-1}. In order to do so, we need to define some new functions and construct a comparison functional below the energy functional $I(v)$.
Recall $g_{1}(s)=\max\{g(s)+ms^q, 0\}$ if $s\ge 0$ and
$g_{1}(s)=-g_{1}(-s)$ if $s\le 0$, and $G_1(s)=\int_{0}^{s}g_1(t)dt$. We set $f(s):=g_3(s)-m|s|^{q-1}s$ and $F(s):=\int_{0}^{s}f(t)dt $ for $s\in \RR$, where $g_3$ is defined by
\begin{equation*}
g_3(s):=s^{p_c}\sup_{0<\tau \le s}\frac{g_1(\tau)}{\tau^{p_c}},\ \textrm{if $s>0$}; \ \ g_3(0):=0;\ \ g_3(s):=-g_3(-s),\ \textrm{if $s<0$},
\end{equation*}
where $p_c \in ( \max(1,q), \ 2^*-1)$ for $N\ge 3 $, and $p_c=q+1$ for $N=2$. Observing that by (g1), there exists a small $\delta>0$ such that $ g_1(s)=G_1(s)=0$ if $|s|\le \delta$, we see that $g_3(s)$ and $f(s)$ are odd continuous functions. Actually $f(s)$ shares similar properties satisfied by $g$.
\begin{lemma}\label{growth-1}
\begin{enumerate}
\item [(i)] $g_3(s) \ge 0$ if $s\ge 0$, $g_3(s)=0$ if $|s|\le \delta$, and $\frac{g_3(s)}{s^{p_c}}$ is increasing for $s >0$.
\item [(ii)] $ f(s)\ge g(s) $ for $s\ge 0$ and $ F(s)\ge G(s) $ for $s\in \RR$. Moreover $\lim_{s \to 0^{+}} \frac{f(s)}{s^q}=-m$ and $F(s_0)>0$, where $s_0$ is from (g3).
\item [(iii)] If $N=2$, for any $\alpha>0$, there exists $C_{\alpha}>0$ such that $|g_3(s)| \le C_{\alpha}e^{\alpha s^2}$ for $s\in \RR$; \\ if $N\ge 3$, $ \lim_{s\to \infty}|g_3(s)| / |s|^{2^{*}-1}=0$.
\item [(iv)] $0\le (p_c+1)G_3(s)\le s g_3(s) $ for $s >0$, where $G_3(s)=\int_{0}^{s}g_3(t)dt $.
\end{enumerate}
\end{lemma}
\begin{proof}
By (g1) and definition of $g_3$, (i) is trivial. From definition of $g_1$ and $g_3$, when $s>0$, we have $g_3(s)\ge g_1(s) \ge g(s)+ms^q$, which shows that $f(s)\ge g(s)$ for $s\ge 0$. From (i) we get $\lim_{s \to 0^{+}} \frac{f(s)}{s^q}=-m$. By (ii) and (g3), we also have $F(s_0)>0$, where $s_0$ is from (g3). So (ii) is also clear.

(iii) Let $N=2$ first. By (g2), for
$\alpha>0$, there a constant $C$ such that $|g(s)|\le C e^{(\alpha/2)s^2}$ for all $s\in \RR$. Therefore $|g_1(s)|\le C_1 e^{(\alpha/2)s^2}$ for all $s\in \RR$. From (i) and definition of (g3), for $s >0$, we have
$$0\le  g_3(s)=s^{p_c}\sup_{0<\tau \le s}\frac{g_1(\tau)}{\tau^{p_c}}\le C_1 \frac{s^{p_c}}{\delta^{p_c}}e^{(\alpha/2)s^2}\le C_2e^{\alpha s^2}.$$
Next let $N \ge 3$. By (g2), we know, for any $\epsilon>0$, there exists $s_\epsilon >0 $ such that
$ g_1(s)/s^{2^*-1} \le \epsilon $ if $ s\ge s_\epsilon$. By continuity, there exists $\delta \le \tau_s \le s$ such that
\begin{equation*}
\frac{g_1(\tau_s)} {\tau_s^{p_c}}=\max_{0<\tau \le s}\frac{ g_1(\tau)} {\tau^{p_c}}.
\end{equation*}
 Set $A_{s_\epsilon}:=\max_{0<\tau \le s_\epsilon} g_1(\tau)$. Now for $s\ge s_\epsilon+(A_{s_\epsilon}/\delta^{p_c}\epsilon)^{1/(2^*-1-p_c)}$, we have
\begin{equation*}
\frac{g_3(s)}{s^{2^{*}-1}}=\frac{s^{p_c}}{s^{2^{*}-1}} \frac{g_1(\tau_s)}{\tau_s^{p_c}} \le \frac{\epsilon }{\delta^{p_c}}\ \ \textrm{if} \ \tau_s \le s_\epsilon \ \textrm{and} \ \ \frac{g_3(s)}{s^{2^{*}-1}}=\frac{\tau_s^{2^*-1-p_c}}{s^{2^{*}-1-p_c}} \frac{g_1(\tau_s)}{\tau_s^{2^{*}-1}} \le \epsilon \ \ \textrm{if} \ \tau_s > s_\epsilon .
\end{equation*}
So we conclude $ \lim_{s\to \infty}|g_3(s)| / |s|^{2^{*}-1}=0$.

(iv) Finally from the definition of $g_3(s)$ and $G_3(s)$, for $s >0$, we obtain
$$ G_3(s)=\int_{0}^{s}g_3(t)dt  \le \sup_{0<\tau \le s}\frac{g_1(\tau)}{\tau^{p_c}}  \int_{0}^{s}t^{p_c} dt=\frac{s}{p_c+1}s^{p_c}\sup_{0<\tau \le s}\frac{g_1(\tau)}{\tau^{p_c}}=\frac{s}{p_c+1}g_3(s),$$
which shows that (iv) is also true.
\end{proof}
Now we define the comparison functional $I_c(v)$ from $E$ to $\RR$ by
\begin{equation*}
I_{c}(v)= \frac{1}{2} \int_{\RR^{N}} |\nabla{v}|^{2}dx- \int_{\RR^{N}}F(v)dx=\frac{1}{2}||\nabla{v}||_2^2+\frac{m}{q+1}||v||_{q+1}^{q+1}-\int_{\RR^{N}}G_3(v)dx, \quad v \in E.
\end{equation*}
It is easy to see that critical points of $I_c(v)$ are solutions of the nonlinear problem $-\triangle{u}=f(u)$ in $\RR^N$ with $u \in E$.
\begin{lemma}\label{mountain-pass-structure-2}
\begin{enumerate}
\item [(i)] $I_c(v) \in C^{1}(E, \RR)$ and $I(v)\ge I_c(v)$ for any $v \in E$.
\item [(ii)] There exist $a_4>0$ and $\rho_4>0$ such that $$I(v)\ge I_c(v)\ge 0, \ \forall \ ||v||_E \le \rho_4 ; \ \ I(v)\ge I_c(v) \ge a_4 , \  \forall \ ||v||_E = \rho_4.$$
\end{enumerate}
\end{lemma}
\begin{proof}
By lemma~\ref{growth-1}, we find out that $f(s)$ satisfy similar version of (g1),(g2) and (g3). Therefore we follow the proof of theorem~\ref{c1functional} and get $I_c(v) \in C^{1}(E, \RR)$, and repeat the proof for
lemma~\ref{mountain-pass-structure-1} (i) and obtain the existence of $a_4$ and $\rho_4$. Since $ F(s)\ge G(s) $ for $s\in \RR$, we have $I(v)\ge I_c(v)$ for any $v \in E$.
\end{proof}
Lemma~\ref{mountain-pass-structure-2} implies $I_c(v)$ has a mountain pass structure around $v=0$. Like we did in section 3, we define a sequence of minimax values $d_n^H$:
\begin{equation}\label{minimax-3}
d_{n}^H:=\inf_{\gamma \in \Phi_{n}^H} \sup_{x\in D_n} I_c(\gamma(x)),
\end{equation}
where $\Phi_n^H:=\{ \gamma\in C(D_n, H)\ |\ \gamma(-x)=-\gamma(x),\ x\in D_n\ ; \ I_c(\gamma(x))<0,\ x\in \partial{D_n} \}. $
\begin{lemma}
For all $n \ge 1$, $\Phi_n^H \neq \emptyset $, $d_{n+1}^H\ge d_n^H$ and $c_n^H \ge d_n^H \ge a_4$, where $a_4$ is from lemma~\ref{mountain-pass-structure-2}.
\end{lemma}
\begin{proof}
Since $f(s)$ satisfy similar version of (g1),(g2) and (g3), repeat the proof of lemma~\ref{nonempty-1} and we have $\Phi_n^H \neq \emptyset $, $d_{n+1}^H\ge d_n^H\ge a_4$.
By lemma~\ref{mountain-pass-structure-2}, $I(v)\ge I_c(v)$ for any $v \in E$,
which implies $c_n^H \ge d_n^H \ge a_4$.
\end{proof}
For each $d_n^H$, according to standard deformation results, we obtain a Palais-Smale sequence. By lemma~\ref{growth-1}(iv), the nonlinear function $g_3$ satisfies Ambrosetti-Rabinowitz condition.
Since the nonlinear structure of $I_c$ is clear, unlike what we did to $c_n^H$ from \eqref{minimax-1} in section 3, we can easily show the compactness of the corresponding Palais-Smale sequence.
More generally, for a sequence of functions $\{ u_i \}\subset H$ such that, as $i \to \infty$,
\begin{equation}\label{Palais-Smale-3}
I_c(u_i)\ \textrm{is uniformly bounded} \ \ \textrm{and} \ \ \ I_{c} '(u_i) \to 0 \ \ \textrm{in $H^*$ strongly},
\end{equation}
we can prove the following lemma.
\begin{lemma}\label{compactness-3}
The Palais-Smale sequence $\{u_i\}$ from \eqref{Palais-Smale-3} is compact in $H$.
\end{lemma}
\begin{proof}
Set $\delta_i = ||I_c'(u_i)||_{H^*}$, by assumption we have $\lim_{i\to \infty}\delta_i=0$, and there exists $C>0$ such that  $|I_c(u_i)|\le C$ uniformly. From lemma~\ref{growth-1}(iv) we find that
\begin{multline*}
C+ \delta_i ||u_i||_{E} \ge I_c(u_i)-\frac{1}{p_c+1}I_c'(u_i)u_i=(\frac{1}{2}-\frac{1}{p_c+1})||\nabla{u_i}||_2^2+(\frac{m}{q+1}-\frac{m}{p_c+1})||u_i||_{q+1}^{q+1} \\
+\int_{\RR^N} \frac{1}{p_c+1}g_3(u_i)u_i-G_3(u_i)dx \ge (\frac{1}{2}-\frac{1}{q+2})||\nabla{u_i}||_2^2+m(\frac{1}{q+1}-\frac{1}{p_c+1})||u_i||_{q+1}^{q+1},
\end{multline*}
which means that $||u_i||_{E}$ is uniformly bounded, so for some $M>0$, we have $||u_i||_{E} \le M$ for all $i$. Taking a subsequence if necessary,
we may assume $u_i \rightharpoonup w$ weakly in $E$ and pointwise a.e. in $\RR^N$ as $i \to \infty$.
Similar to similar to \eqref{zero-3}, \eqref{TM-2} and \eqref{TM-3}, by weak convergence, we have
\begin{equation*}
\int_{\RR^N}f(u_i) vdx \to \int_{\RR^N}f(w) vdx \ \ \textrm{and} \ \ I_c'(w)v =0 \ \textrm{for any $v \in H$},
\end{equation*}
which says that $I_c'(w)w=0$.
Now we can repeat the proof of theorem~\ref{compactness-1} with $g_3(s)$ replacing $g_1(s)$, and we obtain, as $i\to \infty$,
\begin{equation}\label{growth-2}
\int_{\RR^N}g_3(u_i)u_idx \to \int_{\RR^N}g_3(w)wdx.
\end{equation}
With \eqref{growth-2} and  $I_c'(w)w=0$, we get that
\begin{multline*}
||\nabla{u_i}||_2^2+m||u_i||_{q+1}^{q+1}=I_c'(u_i)u_i+\int_{\RR^N}g_3(u_i)u_idx \to  \int_{\RR^N}g_3(w)wdx=||\nabla{w}||_2^2+m||w||_{q+1}^{q+1},
\end{multline*}
which implies the strong convergence of $u_i$ to $w$ in $H$.
\end{proof}
\begin{remark}
With the principle of symmetric criticality \cite{PAL}, the compactness lemma~\ref{compactness-3} guarantees that there exists $v_n^H \in H$ such that $I_c(v_n^H)=d_n^H$ and $I_c'(v_n^H)=0$, but this is not the main point here. In fact we can go through the same process of section 3 to show that $d_n^H$ is achieved by some critical point of $I_c(v)$. By lemma~\ref{growth-1}, there exist $C_1>0$ and $C_2>0$ such that
\begin{equation*}
G_3(s)\ge C_1 |s|^{p_c+1}-C_2|s|^{q+1}, \ \forall s \in \RR.
\end{equation*}
Because all norms are equivalent in finite dimensional space, for any finite dimensional subspace $X \subset H$, there exists $R(X)>0$ such that
\begin{equation}\label{growth-3}
I_c(v)\le \frac{1}{2}||\nabla{v}||_2^2+(\frac{m}{q+1}+C_2)||v||_{q+1}^{q+1}- C_1 ||v||_{q+2}^{p_c+1} <0 , \ \forall v \in X \ \textrm{ with $||v||_E\ge R(X)$}.
\end{equation}
The main point is that, by theorem 9.12 in \cite{RAB}, lemma~\ref{mountain-pass-structure-2} along with lemma~\ref{compactness-3} and \eqref{growth-3} guarantees that $I_c(v)$ possesses an unbounded sequence of critical values.
We are going to repeat the proof of theorem 9.12 in \cite{RAB} under our framework, and we only outline the main steps through lemma~\ref{unbounded-1} and theorem~\ref{unbounded-2}. For more details, check \cite{RAB}.
\end{remark}
We define a sequence of minimax values $e_n^H$:
\begin{equation}\label{minimax-4}
e_{n}^H:=\inf_{B \in \Psi_{n}^H} \sup_{v\in B} I_c(v),
\end{equation}
where $\Psi_n^H:=\cup_{l=n}^{\infty} \{ \gamma(\overline{D_l - A})\ | \ \gamma \in \Phi_l^H; \ A\in \mathcal{A}_l \ ; \ gen(A)\le l-n \} $,
$\mathcal{A}_l$ is the family of closed sets $A \subset \RR^{l}-\{0\}$ such that $A=-A$ and $gen(A)$ is the Krasnoselskii genus of $A$. Noticing that $A$ is empty if $gen(A)=0$, so we have the following lemma.
\begin{lemma}\label{unbounded-1}
Let $n$ and $i$ be two positive integers. Assume that $Z \subset H-\{0\}$ is a closed set such that $Z=-Z$ and $gen(Z)<\infty$.
\begin{enumerate}
\item [(i)] For all $n\ge 1$, $\{\gamma(D_n) \ |\ \gamma \in \Phi_n^H \} \subset \Psi_n^H$ and $\Psi_{n+1}^H \subset \Psi_n^H$. In particular $d_n^H\ge e_n^H$ and $e_{n+1}^H \ge e_n^H$.
\item [(ii)] For all $n\ge 1$, if $B\in \Psi_n^H$, then $B\cap \{ v\in E \ |\ ||v||_E =\rho_4 \} \neq \emptyset $ and $e_n^H \ge a_4$, here $a_4$ and $\rho_4$ are from lemma~\ref{mountain-pass-structure-2}.
\item [(iii)] If $gen(Z)\le i <n$ and $B \in \Psi_n^H$, then $\overline{B-Z} \in \Psi_{n-i}^H$.
\item [(iv)] If $e_{n}^H=\cdots=e_{n+i}^H$, then $gen(K_{e_n^H})\ge i+1$, where the symmetric set $K_{e_n^H}=\{ v\in H \ | \ I_c(v)=e_n^H; \ I_c'(v)=0 \} $.
\end{enumerate}
\end{lemma}
\begin{theorem}\label{unbounded-2}
Assume that (g1),(g2) and (g3) hold. There exists a sequence of critical points $\{v_n^H\}\subset H$ of $I_c(v)$ with $I_c(v_n^H)=e_n^H$ and $ \lim_{n\to \infty}e_n^H = \infty$. In particular $\lim_{n\to \infty}c_n^H = \infty$.
\end{theorem}
\begin{proof}
Since $e_n^H\ge a_4 >0$, we can construct an odd continuous pseudo-gradient flow $\eta(\tau,v)$ such that $\eta(1,v)=v$ for any $v$ with $I_c(v)<0$. So we have $\eta(1,B)\in \Psi_n^H$ for any $B \in \Psi_n^H$, which implies the existence of Palais-Smale sequence at each minimax value $e_n^H$. By lemma~\ref{compactness-3}, $e_n^H$ is achieved by $I_c(v)$ at $v_n^H$. Since $c_n^H \ge d_n^H \ge e_n^H$, we are left to show $ \lim_{n\to \infty}e_n^H = \infty$.
Assume otherwise $e_n^H$ is uniformly bounded. Noticing that $e_n^H$ is increasing, then we set $e^*=\lim_{n\to \infty}e_n^H$. If $e_n^H =e^*$ for large $n$, lemma~\ref{unbounded-1}(iv) implies that $gen(K_{e^*})=\infty$.
But by lemma~\ref{compactness-3} $K_{e^*}$ is compact, so $gen(K_{e^*})$ is finite. Therefore we must have $e^{*} >e_n^H$ for any $n\ge 1$. Let $K^*=\{ v\in H \ | \ e_1\le I_c(v)\le e^*; \ I_c'(v)=0 \} $,
again by lemma~\ref{compactness-3} $K^{*}$ is compact and $gen(K^{*})$ is finite. Set $i^*=gen(K^*)$. There is a small $\delta>0$ such that  $gen(K^*_\delta)=i^*$, here
$K^*_\delta=\{ v\in H \ | \ \inf_{w\in K^*} ||v-w||_E \le \delta \} $. For small $\epsilon>0$, there exists $n^*>0$ such that $e_n^H\in (e^*-\epsilon/2, e^*)$ if $n \ge n^*$. By definition of $e_{n^*+i^*}^H$,
there is $B\in \Psi_{n^*+i^*}^H$ such that $\max_{v \in B} I_c(v) \le e_{n^*+i^*}^H+\epsilon$. It is clear that $\max_{v \in B} I_c(v) \le e^*+\epsilon$ and $\overline{B-K^*_\delta} \in \Psi_{n^*}^H$.
We can construct an odd continuous pseudo-gradient flow $\eta(\tau,v)$ such that $I_c(\eta(1,\overline{B-K^*_\delta}))\le e^*-\epsilon$ and $\eta(1,v)=v$ for any $v$ with $I_c(v)<0$.
Since $\eta(1,\overline{B-K^*_\delta}) \in \Psi_{n^*}^H$, we have
$$ e_{n^*}^H \le \max_{v \in \overline{B-K^*_\delta}} I_c(\eta(1, v)) \le e^*-\epsilon \le e_{n^*}^H+\epsilon/2 -\epsilon = e_{n^*}^H-\epsilon/2,$$
which is a contraction. So we have  $ \lim_{n\to \infty}e_n^H = \infty$.
\end{proof}



\end{document}